\documentclass[reqno,A4paper]{amsart}

\usepackage{amsmath}
\usepackage{amssymb}
\usepackage{amsthm}
\usepackage{enumerate}
\usepackage{mathrsfs} 
\usepackage{eqlist}
\usepackage{array}

\usepackage{tikz-cd}
\usepackage{acadpgf}

\theoremstyle{plain}
\newtheorem{theorem}{Theorem}[section]

\newtheorem{lemma}{Lemma}[section]
\newtheorem{corol}{Corollary}[theorem]

\newtheorem*{OPEN}{Open problem}

\theoremstyle{definition}
\newtheorem{definition}{Definition}
\newtheorem{example}{\textit{Example}} 

\numberwithin{equation}{section}

\newcommand{\NN}{\mathbb{N}}    

\newcommand\AAA{\mathbf{A}}
\newcommand\BB{\mathbf{B}}
\newcommand\CCC{\mathbf{C}}
\newcommand\DD{\mathbf{D}}
\newcommand\EE{\mathbf{E}}
\newcommand\EM{\mathbf{EM}}
\newcommand\KK{\mathbf{K}}

\newcommand\calA{\mathcal{A}}
\newcommand\calB{\mathcal{B}}
\newcommand\calC{\mathcal{C}}

\newcommand\calE{\mathcal{E}}

\newcommand\calU{\mathcal{U}}
\newcommand\calV{\mathcal{V}}

\newcommand\0{\varnothing}
\renewcommand\epsilon{\varepsilon}
\renewcommand\phi{\varphi}
\let\le\leqslant
\let\ge\geqslant
\newcommand{\lt}[1]{\mathrel{<^{#1}}}

\newcommand{\lex}[1]{\mathrel{<^{#1}_{\mathit{lex}}}}

\newcommand{\restr}[2]{\mbox{$#1$}\mbox{$\upharpoonright$}_{#2}}

\newcommand\Ob{\mathrm{Ob}}

\newcommand\ID{\mathrm{ID}}
\newcommand\id{\mathrm{id}}

\newcommand{\Boxed}[1]{\hbox{$#1$}}

\newcommand{\Set}{\mathbf{Set}}
\newcommand{\SetEmb}{\mathbf{Set}_{\emb}}
\newcommand{\SetInj}{\mathbf{Set}_{\inj}}
\newcommand{\SetEmbFin}{\mathbf{Set}_{\emb}^{\fin}}
\newcommand{\OsetEmb}{\mathbf{Oset}_{\emb}}
\newcommand{\OsetEmbFin}{\mathbf{Oset}_{\emb}^{\fin}}

\newcommand{\ChEmb}{\mathbf{Ch}_{\emb}}
\newcommand{\ChEmbFin}{\mathbf{Ch}_{\emb}^{\fin}}

\newcommand{\OalgEmb}{\mathbf{Oalg}_{\emb}}
\newcommand{\OalgEmbFin}{\mathbf{Oalg}_{\emb}^{\fin}}

\newcommand{\AlgEmb}{\mathbf{Alg}_{\emb}}
\newcommand{\AlgEmbFin}{\mathbf{Alg}_{\emb}^{\fin}}

\newcommand{\OrfEmb}{\mathbf{Orf}_{\emb}}
\newcommand{\OrfEmbFin}{\mathbf{Orf}_{\emb}^{\fin}}

\newcommand{\RfEmb}{\mathbf{Rf}_{\emb}}
\newcommand{\RfEmbFin}{\mathbf{Rf}_{\emb}^{\fin}}

\newcommand\op{\mathit{op}}

\newcommand\fin{\mathit{fin}}

\newcommand\emb{\mathit{emb}}
\newcommand\inj{\mathit{inj}}

\newcommand\Aut{\mathrm{Aut}}
\newcommand\iso{\mathrm{iso}}

\begin{document}

\title{Coalgebraic methods for Ramsey degrees of unary algebras}
\author{Dragan Ma\v sulovi\'c}
\address{Department of Mathematics and Informatics\\
         Faculty of Sciences, University of Novi Sad, Serbia}
\email{dragan.masulovic@dmi.uns.ac.rs}

\thanks{This research was supported by the Science Fund of the Republic of Serbia, Grant No.~7750027:
Set-theoretic, model-theoretic and Ramsey-theoretic phenomena in mathematical structures: similarity and diversity -- SMART}

\subjclass[2020]{Primary 08A60; Secondary 18C20} 
\keywords{Ramsey property, Ramsey degrees, Eilenberg-Moore category, comonad, unary algebra}

\begin{abstract}
  In this paper we prove the existence of small and big Ramsey degrees of classes of finite
  unary algebras in an arbitrary (not necessarily finite) algebraic language $\Omega$. Our
  results generalize some Ramsey-type results of M.\ Soki\'c concerning finite unary algebras over
  finite languages. To do so we develop a completely new strategy that relies on the fact that right
  adjoints preserve the Ramsey property.
  We then treat unary algebras as Eilenberg-Moore coalgebras for a functor with comultiplication
  and using pre-adjunctions transport the Ramsey properties we are interested in from the category of finite or
  countably infinite chains of order type $\omega$. Moreover, we show that finite objects have finite
  big Ramsey degrees in the corresponding cofree structures over countably many generators.
\end{abstract}

\maketitle

\section{Introduction}

Almost any reasonable class of finite relational structures has some kind of Ramsey property.
In the realm of finite algebras the picture is considerably different -- in the past 50 years
the Ramsey property has been established in only a handful of cases:
finite Boolean algebras~\cite{graham-rothschild-DRT}; finite vector spaces over a
fixed finite field~\cite{graham-leeb-rothschild}; finite Boolean lattices~\cite{promel-85};
finite unary algebras over a finite language~\cite{sokic-unary-functions};
finite $G$-sets for a finite group $G$~\cite{sokic-unary-functions};
and finite semilattices~\cite{sokic-semilattices}. It is notable that,
despite a long history~\cite{deuber-rothschild,oates-williams,voigt-fin-ab-grps}, the following question
is still open:

\begin{OPEN}
  Is it true that the class of finite groups has small Ramsey degrees?
\end{OPEN}

In this paper we generalize some of the results of Soki\'c from~\cite{sokic-unary-functions} and show that
for an arbitrary monoid $M$ (finite or infinite) the class of all finite $M$-sets has finite small Ramsey degrees.
This immediately implies that the class of all finite $G$-sets has finite small Ramsey degrees for an arbitrary group $G$ (finite or infinite),
and that the class of all finite unary algebras over an arbitrary (finite or infinite) algebraic language $\Omega$
has finite small Ramsey degrees. Moreover, we show that in all of the three contexts finite objects have finite
big Ramsey degrees in the corresponding cofree structure over countably many generators.

All our results are spelled out using the categorical reinterpretation of the Ramsey property as proposed in~\cite{masul-scow}.
Actually, it was Leeb who pointed out already in 1970 that the use of category theory can be quite helpful
both in the formulation and in the proofs of results pertaining to structural Ramsey theory~\cite{leeb-cat}.

In order to prove the existence of small and big Ramsey degrees in classes of $M$-sets (and, thus, immediately
get the corresponding results about $G$-sets for an arbitrary group $G$, and $\Omega$-algebras for an arbitrary
unary algebraic language $\Omega$) we use a completely new strategy developed in~\cite{masul-RPEMKLEI-MONAD}
that relies on transporting the Ramsey property from a category to the associated category of
Eilenberg-Moore (co)algebras for a (co)monad. In~\cite{masul-RPEMKLEI-MONAD} we show how \emph{dual} Ramsey phenomena
transport from a category to the category of Eilenberg-Moore algebras for a monad.
In this paper we are interested in analogous results for comonads.
Although comonads are categorical duals of monads, the combinatorics involved is of a different kind,
and modeling algebraic phenomena by coalgebras for a comonad is far more challenging.
For a monoid $M$ we consider $M$-sets as Eilenberg-Moore coalgebras for monoid action comonad (to be described below).
The results then follow from the Finite Ramsey Theorem and the Infinite Ramsey Theorem
by building an appropriate endofunctor and a comultiplication on the category of
finite and countably infinite chains of order type~$\omega$.

After fixing standard notions and notation in Section~\ref{rpemklei.sec.prelim}, we present our proof strategy in
Section~\ref{rpemklei.sec.RP}. Our starting point is the observation from~\cite{masul-scow} that right adjoints preserve the
Ramsey property. We then observe that if $E$ is a comonad on a
category with the Ramsey property, both the Kleisli category and the Eilenberg-Moore category for the comonad have the
Ramsey property. Unfortunately, these two simple results are not very useful: for the categorical treatment of the Ramsey property
it is essential to restrict the attention to categories where the morphisms are mono, and counits for comonads
we are interested in cannot be expected to consist of monos for cardinality reasons
(see Example~\ref{rpemklei.ex.list-comonad}). This makes it impossible to model phenomena we are interested in
by comonads over categories whose morphisms are mono. Therefore, we relax the context by proving that
the Ramsey property carries over from the category of chains to the category of \emph{weak Eilenberg-Moore coalgebras}
defined for functors with comultiplication and no counit, which are straightforward weakenings of comonads.
Our main results are presented in Section~\ref{rpemklei.sec.unary-algs}.

\section{Preliminaries}
\label{rpemklei.sec.prelim}

\subsection{Chains}

  A \emph{chain} is a linearly ordered set $(A, \Boxed<)$.
  Finite or countably infinite chains will sometimes be denoted as
  $\{a_1 < \dots < a_n < \dots\}$. For example, $\omega = \{0 < 1 < 2 < \dots\}$.
  Every strict linear order $<$ induces the reflexive version $\le$ in the obvious way.
  A map $f : A \to B$ from a chain $(A, \Boxed<)$ to a chain $(B, \Boxed<)$ is
  \emph{monotone} if $a_1 < a_1 \Rightarrow f(a_1) < f(a_2)$ for all $a_1, a_2 \in A$.

  Let $(A, \Boxed<)$ be a chain. Sequences of elements of $A$ can be \emph{lexicographically ordered} as follows.
  Let $\alpha$ and $\beta$ be ordinals. For a pair of distinct sequences $(a_\xi)_{\xi < \alpha} \in A^\alpha$ and
  $(b_\xi)_{\xi < \beta} \in A^\beta$ (not necessarily of the same length) we write
  $(a_\xi)_{\xi < \alpha} \lex{} (b_\xi)_{\xi < \beta}$ if
  \begin{itemize}
    \item $\alpha \le \beta$ and $a_\xi = b_\xi$ for all $\xi < \alpha$; or
    \item there is a $\xi$ such that $a_\xi \ne b_\xi$ in which case we let $\eta = \min \{\xi : a_\xi \ne b_\xi\}$ and
          declare $(a_\xi)_{\xi < \alpha} \lex{} (b_\xi)_{\xi < \alpha}$ if and only if $a_\eta < b_\eta$.
  \end{itemize}
  In particular, for every chain $(A, \Boxed{<})$ and every well-ordered chain $(S, \Boxed<)$ the set $A^S$ 
  of all functions $f : S \to A$ can be ordered lexicographically and the corresponding lexicographic ordering on $A^S$ will be denoted by $\lex A$.
  For every finite ordinal $n$ and every well-ordered chain $(A, \Boxed{<})$ the chain $(A^n, \Boxed{\lex{A}})$
  is also well-ordered.

\subsection{Unary algebras, $G$-sets and $M$-sets}

  Let $\Omega$ be an algebraic language, that is, the set of constant and functional symbols.
  A \emph{$\Omega$-algebra} is a structure $(A, \Omega^A)$ where $\Omega^A = \{f^A : f \in \Omega\}$ is a set of operations
  on $A$ such that the arity of each operation $f^A$ coincides with
  the arity of the corresponding functional symbol $f \in \Omega$. We say that $A$ is a \emph{carrier} of an
  $\Omega$-algebra $(A, \Omega^A)$.

  A language $\Omega = \{f_i : i \in I \}$ is \emph{unary} if every symbol
  in $\Omega$ is a unary functional symbol; and it is \emph{monounary} if $\Omega$ consists of a single
  unary functional symbol. Correspondingly, a \emph{(mono)unary algebra} is an $\Omega$-algebra where
  $\Omega$ is a (mono)unary algebraic language.

  Let $(M, \Boxed\cdot, 1)$ be a monoid. An \emph{$M$-set} is a pair $(A, \alpha)$
  where $\alpha : M \times A \to A$ is the \emph{structure map} and has the following properties: $\alpha(1, a) = a$ and
  $\alpha(m_1, \alpha(m_2, a)) = \alpha(m_2 \cdot m_1, a)$ for all $a \in A$ and $m_1, m_2 \in M$.
  A mapping $f : A \to B$ is a morphism between $M$-sets $(A, \alpha)$ and $(B, \beta)$ if
  $f(\alpha(m, a)) = \beta(m, f(a))$ for all $a \in A$ and all $m \in M$. An injective
  morphism of $M$-sets will be referred to as an \emph{embedding}. A \emph{cofree $M$-set on the
  set of generators~$X$} is the $M$-set $\calE(X) = (X^M, \gamma)$ where $\gamma(m, h) \in X^M$ is
  given by $\gamma(m, h)(m') = h(m \cdot m')$. See Lemma~\ref{rpemklei.lem.uniq-cofree} and
  Example~\ref{rpemklei.ex.mon-act} below. In particular, a \emph{$G$-set} is an $M$-set
  where $M = G$ is a group.
  
  Let $\Omega = \{f_i : i \in I \}$ be an arbitrary unary algebraic language.
  Let $M = \Omega^*$ be the free monoid of finite words over $\Omega$ and let $1 \in \Omega^*$ denote the empty word.
  Every $\Omega$-algebra $\calA = (A, \Omega^A)$ can be represented by an $M$-set
  $(A, \alpha)$ where the structure map $\alpha : M \times A \to A$ evaluates words from $M = \Omega^*$ in $\calA$ in the obvious way:
  $\alpha(1, a) = a^1 = a$ and $\alpha(f_{i_1} \dots f_{i_k}, a) = a^{f_{i_1} \dots f_{i_k}} = (f^A_{i_k} \circ \dots \circ f^A_{i_1})(a)$.
  Note that a map is a homomorphism between two unary algebras if and only if it is a morphism of corresponding $M$-sets.

\subsection{Categories and functors}

For basic category-theoretic notions and notation (category, functor, homset, natural transformation etc)
we refer the reader to~\cite{AHS}. Let us here quickly fix some notation.
Let $\CCC$ be a category. By $\Ob(\CCC)$ we denote the class of all the objects in $\CCC$.
Homsets in $\CCC$ will be denoted by $\hom_\CCC(A, B)$, or simply $\hom(A, B)$ when $\CCC$ is clear from the context.
The identity morphism will be denoted by $\id_A$ and the composition of morphisms by $\Boxed\cdot$ (dot).
If $\hom_\CCC(A, B) \ne \0$ we write $A \to B$.

A morphism $f \in \hom_\CCC(A, B)$ is \emph{mono} if it is left cancellable in the following sense:
for all $C \in \Ob(\CCC)$ and all $g, h \in \hom_\CCC(C, A)$, if $f \cdot g = f \cdot h$ then $g = h$;
and it is \emph{iso} or \emph{invertible} if there exists an $f' \in \hom_\CCC(B, A)$ such that 
$f \cdot f' = \id_B$ and $f' \cdot f = \id_A$.
Let $\iso_\CCC(A, B)$ denote the set of all the invertible morphisms $A \to B$, and let
$\Aut_\CCC(A) = \iso(A, A)$ denote the set of all the invertible morphisms $A \to A$.
An object $A \in \Ob(\CCC)$ is \emph{rigid} if $\Aut_\CCC(A) = \{\id_A\}$.
As usual, $\CCC^\op$ denotes the opposite category.
If $\CCC$ is a category of structures, where by a structure we mean a set together with some additional information, by
$\CCC^{\fin}$ we denote the full subcategory of $\CCC$ spanned by its finite members.

A category $\DD$ is a \emph{subcategory} of a category $\CCC$ if $\Ob(\DD) \subseteq\Ob(\CCC)$ and
$\hom_\DD(A, B) \subseteq \hom_\CCC(A, B)$ for all $A, B \in \Ob(\DD)$; and it is a \emph{full subcategory}
if it a subcategory of $\CCC$ such that $\hom_\DD(A, B) = \hom_\CCC(A, B)$ for all $A, B \in \Ob(\DD)$.
Let $\DD$ be a full subcategory of $\CCC$. An $S \in \Ob(\CCC)$ is \emph{universal for $\DD$} if for every $D \in \Ob(\DD)$
the set $\hom_\CCC(D, S)$ is nonempty and consists of monos only.
Note that if there exists an $S \in \Ob(\CCC)$ universal for $\DD$ then all the morphisms
in $\DD$ are mono. If $\DD$ is a full subcategory of $\CCC$ the we say that $\CCC$ is an \emph{ambient category} for $\DD$.
An ambient category $\CCC$ is usually a category in which we can perform certain operations
that are not possible in $\DD$, or which contains an object universal for~$\DD$.

A \emph{functor} $F : \CCC \to \DD$ from a category $\CCC$ to a category $\DD$ maps $\Ob(\CCC)$ to
$\Ob(\DD)$ and maps morphisms of $\CCC$ to morphisms of $\DD$ so that
$F(f) \in \hom_\DD(F(A), F(B))$ whenever $f \in \hom_\CCC(A, B)$, $F(f \cdot g) = F(f) \cdot F(g)$ whenever
$f \cdot g$ is defined, and $F(\id_A) = \id_{F(A)}$.
A functor $F : \CCC \to \DD$ is \emph{faithful} if it is injective on homsets in the sense that
for all $A, B \in \Ob(\CCC)$ and all $f, g \in \hom_\CCC(A, B)$, if $F(f) = F(g)$ then $f = g$.
A functor $U : \CCC \to \DD$ is \emph{forgetful} if it is faithful and surjective on objects
(that is, for every $D \in \Ob(\DD)$ there is a $C \in \Ob(\CCC)$ with $U(C) = D$).

If $U : \CCC \to \DD$ is a forgetful functor
we may actually assume that $\hom_{\CCC}(A, B) \subseteq \hom_\DD(U(A), U(B))$ for all $A, B \in \Ob(\CCC)$.
The intuition behind this point of view is that $\CCC$ is a category of objects with a lot of structure,
say topological groups, $\DD$ is a category of objects with less structure, say groups,
and $U$ takes a ``rich'' object $\calA$ and ``forgets'' the extra structure (topology, say) to produce
its ``poor'' version~$A$.
Then for every morphism $f : \calA \to \calB$ in $\CCC$ the same map is a morphism $f : A \to B$ in $\DD$.
Therefore, if $U$ is a forgetful functor we shall always take that $U(f) = f$. In particular,
$U(\id_{A}) = \id_{U(A)}$ and we, therefore, identify $\id_{A}$ with $\id_{U(A)}$.
Also, if $U : \CCC \to \DD$ is a forgetful functor and all the morphisms in $\DD$ are mono, then
all the morphisms in $\CCC$ are mono.

\begin{example}
  \begin{enumerate}[(1)]
  \item
    $\Set$ is the category of sets and all set functions. $\SetInj$ is the category of sets and injective functions;
    this is a subcategory of $\Set$, although not a full one.

  \item
    $\ChEmb$ is the category whose objects are chains and whose morphisms are embeddings.

  \item
    For a monoid $M$ let $\SetEmb(M)$ denote the category of $M$-sets and embeddings of $M$-sets,
    and for a group $G$ let $\SetEmb(G)$ denote the category of $G$-sets and embeddings.
  
  \item
    For an algebraic language $\Omega$ let $\AlgEmb(\Omega)$ denote the category of $\Omega$-algebras and embeddings, and let
    $\OalgEmb(\Omega)$ denote the category of ordered $\Omega$-algebras and embeddings.
  \end{enumerate}
\end{example}

\subsection{Adjunctions, monads and comonads}

An \emph{adjunction} between categories $\BB$ and $\CCC$ consists of a pair of functors $F : \BB \rightleftarrows \CCC : H$
together with a family of isomorphisms
$
  \Phi_{X,Y} : \hom_\CCC(F(X), Y) \overset\cong\longrightarrow \hom_\BB(X, H(Y))
$
indexed by pairs $(X, Y) \in \Ob(\BB) \times \Ob(\CCC)$ and
natural in both $X$ and $Y$. The functor $F$ is then \emph{left adjoint} (to $H$) and $H$ is \emph{right adjoint} (to $F$).

Let $\CCC$ be a category and $T : \CCC \to \CCC$ a functor. By dualizing the notions of multiplication for an
endofunctor and the monad we arrive at the notions of the comultiplication and comonad as follows.
A \emph{comultiplication for an endofunctor} $E : \CCC \to \CCC$ is a natural transformation
$\delta : E \to EE$ such that for each $A \in \Ob(\CCC)$ the diagram on the left commutes:
\begin{center}
    \begin{tikzcd}
      E(A) \arrow[r, "\delta_A"] \arrow[d, "\delta_A"'] & EE(A) \arrow[d, "\delta_{E(A)}"] \\
      EE(A) \arrow[r, "E(\delta_A)"'] & EEE(A)
    \end{tikzcd}
    \quad
    \begin{tikzcd}
      E(A)  & EE(A) \arrow[l, "E(\epsilon_A)"'] \arrow[r, "\epsilon_{E(A)}"] & E(A) \\
            & E(A) \arrow[ul, "\id_{E(A)}"] \arrow[ur, "\id_{E(A)}"'] \arrow[u, "\delta_A"']
    \end{tikzcd}
\end{center}
A natural transformation $\epsilon : E \to \ID$ is a \emph{counit for $\delta$} if the diagram on the right above
commutes for each $A \in \Ob(\CCC)$. A \emph{comonad} on a category $\CCC$ is a triple $(E, \delta, \epsilon)$ where $E : \CCC \to \CCC$ is a functor,
$\delta$ is a comultiplication for $E$ and $\epsilon$ is a counit for~$\delta$.

Let $F : \CCC \to \CCC$ be a functor. An \emph{$F$-algebra} is a pair $(A, \alpha)$ where $\alpha \in \hom_\CCC(F(A), A)$, while
an \emph{$F$-coalgebra} is a pair $(A, \alpha)$ where $\alpha \in \hom_\CCC(A, F(A))$.
An \emph{algebraic homomorphism} between $F$-algebras $(A, \alpha)$ and $(B, \beta)$ is a morphism $f \in \hom_\CCC(A, B)$ such that
the diagram on the left commutes:
\begin{center}
  \begin{tabular}{c@{\qquad}c}
    \begin{tikzcd}
      F(A) \arrow[r, "F(f)"] \arrow[d, "\alpha"'] & F(B) \arrow[d, "\beta"]\\
      A \arrow[r, "f"'] & B
    \end{tikzcd}
    &
    \begin{tikzcd}
      A \arrow[r, "f"] \arrow[d, "\alpha"'] & B \arrow[d, "\beta"]\\
      F(A) \arrow[r, "F(f)"'] & F(B)
    \end{tikzcd}
  \end{tabular}
\end{center}
A \emph{coalgebraic homomorphism} between $F$-coalgebras $(A, \alpha)$ and $(B, \beta)$ is a morphism $f \in \hom_\CCC(A, B)$ such that
the diagram on the right commutes.

To each comonad $(E, \delta, \epsilon)$ we can straightforwardly assign the Kleisli category $\KK = \KK(E, \delta, \epsilon)$
and the Eilenberg-Moore category $\EM = \EM(E, \delta, \epsilon)$ by dualizing the corresponding standard constructions
for a monad. Explicitly,
the objects of the Kleisli category $\KK(E, \delta, \epsilon)$ are the same as the objects of $\CCC$, morphisms are defined by
$
  \hom_\KK(A, B) = \hom_\CCC(E(A), B)
$
and the composition in $\KK$ for $f \in \hom_\KK(A, B)$ and $g \in \hom_\KK(B, C)$ is given by
$
  g \mathbin{\cdot_\KK} f = g \cdot E(f) \cdot \delta_A
$.
The objects of the Eilenberg-Moore category $\EM(E, \delta, \epsilon)$ are Eilenberg-Moore $E$-coalgebras
(special $E$-coalgebras to be defined immediately), morphisms are coalgebraic homomorphisms and the composition is as in $\CCC$.
An \emph{Eilenberg-Moore $E$-coalgebra} is an $E$-coalgebra for which the following two diagrams commute:
\begin{center}
  \begin{tabular}{c@{\qquad}c}
    \begin{tikzcd}
      A \arrow[r, "\alpha"] \arrow[d, "\alpha"']  & E(A) \arrow [d, "\delta_A"]\\
      E(A) \arrow[r, "E(\alpha)"'] & EE(A)
    \end{tikzcd}
    &
    \begin{tikzcd}
      A & E(A) \arrow[l, "\epsilon_A"'] \\
      & A \arrow[u, "\alpha"'] \arrow[ul, "\id_A"]
    \end{tikzcd}
  \end{tabular}
\end{center}
An \emph{weak Eilenberg-Moore $E$-coalgebra} is an $E$-coalgebra for which only the diagram on the left commutes.
Let $\EM^w(E, \delta)$ denote the category of weak Eilenberg-Moore $E$-coalgebras and coalgebraic homomorphisms.

\begin{example}\label{rpemklei.ex.list-comonad}
  The \emph{list comonad} is a comonad that provides the necessary infrastructure for dealing with finite lists.
  For any set $A$ let $E(A) = A^+$ be the set of all the nonempty finite sequences of elements of~$A$,
  and for a set-function $f : A \to B$ let $E(f)$ act coordinatewise: $E(f)(a_1, \dots, a_n) = (f(a_1), \dots, f(a_n))$.
  Define $\delta_A : E(A) \to EE(A)$ by
  $
   \delta_A(a_1, \dots, a_n) = \big(
     (a_1, a_2, \dots, a_n), (a_2, \dots, a_n), \dots, (a_n)
   \big)
  $
  and $\epsilon_A : E(A) \to A$ by $\epsilon_A(a_1, \dots, a_n) = a_1$. Then it is easy to check that $(E, \delta, \epsilon)$
  is a $\Set$-comonad. The Eilenberg-Moore coalgebras for this comonad correspond to rooted forests
  (see Section~\ref{rpemklei.sec.RP}). Note that $\epsilon_A$ cannot be mono for cardinality reasons.
  So, in order to model algebraic phenomena we have to work with functors with comultiplication and no counit.
\end{example}

\begin{example}\label{rpemklei.ex.mon-act}
  Fix a monoid $M$. Each $M$-set can be represented as an algebra for a monad as well as a coalgebra for a comonad.
  Define $T : \Set \to \Set$ by $T(A) = M \times A$ on objects,
  while for a mapping $f : A \to B$ let $T(f) : M \times A \to M \times B
  : (m, a) \mapsto (m, f(a))$. Define $\mu_A : TT(A) \to T(A)$ by $\mu_A(m_1, (m_2, a)) = (m_2 \cdot m_1, a)$ and $\eta_A : A \to T(A)$ by
  $\eta(a) = (1, a)$. Then $(T, \mu, \eta)$ is a monad whose Eilenberg-Moore algebras correspond precisely to $M$-sets.
  Using the natural isomorphism $\hom(M \times A, B) \cong \hom(A, B^M)$,
  $M$-sets can be represented by Eilenberg-Moore coalgebras for the following comonad.
  Define $E : \Set \to \Set$ by $E(A) = A^M$ on objects, while for a mapping $f : A \to B$
  let $E(f) : A^M \to B^M : h \mapsto f \circ h$.
  Next, define $\delta_A : E(A) \to EE(A)$ by $\delta_A(h)(m_1)(m_2) = h(m_2 \cdot m_1)$
  and let $\epsilon_A : E(A) \to A : h \mapsto h(1)$.
  Then $(E, \delta, \epsilon)$ is a comonad whose Eilenberg-Moore coalgebras correspond precisely to $M$-sets.
  Finally, recall that for a unary algebraic language $\Omega$ 
  every $\Omega$-algebra $\calA = (A, \Omega^A)$ can be represented by an $M$-set $\alpha : M \times A \to A$
  where $M = \Omega^*$ is the free monoid of finite words over $\Omega$, and hence by an $E$-coalgebra for the comonad
  we have just described. Then it easily follows that $\AlgEmb(\Omega) \cong \EM(E, \delta, \epsilon)$.
\end{example}

A \emph{cofree Eilenberg-Moore $E$-coalgebra over a set of generators $X$} is the
Eilenberg-Moore $E$-coalgebra $\calE(X) = (E(X), \delta_X)$. The following lemma
motivates the choice of the terminology. We will not use this lemma later in the paper -- the reason we include it is just
to motivate the name.

\begin{lemma}\label{rpemklei.lem.uniq-cofree}
  Let $\calA = (A, \alpha)$ be an Eilenberg-Moore $E$-coalgebra and $X$ a set.
  For every mapping $f : A \to X$ there is a unique coalgebra homomorphism
  $f^\# : \calA \to \calE(X)$ such that $\epsilon_X \cdot f^\# = f$.\qed
\end{lemma}

\begin{example}
  Let us describe the cofree coalgebra on $\omega$ generators for the monoid action comonad and, in particular,
  the cofree unary coalgebra on $\omega$ generators for a unary language $\Omega$.
  Given a monoid $M$, the carrier of the cofree coalgebra on $\omega$ generators is the set $C = \omega^M$ and $M$ acts on $C$
  so that $(m \cdot f)(x) = f(xm)$ for all $f \in C$ and $x, m \in M$. Now, given a unary algebraic language $\Omega$
  the carrier of the cofree $\Omega$-coalgebra on $\omega$ generators is $D = \omega^{\Omega^*}$ and for an $s \in \Omega$
  the corresponding unary operation is given by $s(f)(w) = f(ws)$ for all $f \in D$ and $w \in \Omega^*$.
\end{example}

\section{Ramsey properties in a category}
\label{rpemklei.sec.RP}

In this section we collect and prove several results about the Ramsey property, Ramsey degrees,
dual Ramsey property and dual small Ramsey degrees in a category. We then use the results of this section
as the main tool to obtain new Ramsey results about $G$-sets in Section~\ref{rpemklei.sec.unary-algs}.

For $k \in \NN$, a $k$-coloring of a set $S$ is any mapping $\chi : S \to k$, where,
as usual, we identify $k$ with $\{0, 1,\dots, k-1\}$.
Let $\CCC$ be a locally small category.
For integers $k \ge 2$ and $t \ge 1$, and objects $A, B, C \in \Ob(\CCC)$ we write
$
  C \longrightarrow (B)^{A}_{k, t}
$
to denote that for every $k$-coloring $\chi : \hom(A, C) \to k$ there is a morphism
$w \in \hom(B, C)$ such that $|\chi(w \cdot \hom(A, B))| \le t$.
(For a set of morphisms $F$ we let $w \cdot F = \{ w \cdot f : f \in F \}$.)
In case $t = 1$ we write
$
  C \longrightarrow (B)^{A}_{k}.
$
We write $C \longleftarrow (B)^{A}_{k, t}$, resp.\ $C \longleftarrow (B)^{A}_{k}$,
to denote that $C \longrightarrow (B)^{A}_{k, t}$, resp.\ $C \longrightarrow (B)^{A}_{k}$, in $\CCC^\op$.

\begin{lemma}\label{rpemklei.lem.C-D} \cite[Lemma 2.4]{masul-scow}
  Let $\CCC$ be a locally small category such that all the morphisms in $\CCC$ are mono and let $A, B, C, D \in \Ob(\CCC)$. If
  $C \longrightarrow (B)^A_{k,t}$ for some $k, t \ge 2$ and if $C \overset\CCC\longrightarrow D$, then $D \longrightarrow (B)^A_{k,t}$.\qed
\end{lemma}

The above lemma tells us that with the arrow notation we can always go ``up to a superstructure'' of $C$.
In some cases we can also ``go down to a substructure'' of $C$.
Let $\AAA$ be a subcategory of $\CCC$ and let $C \in \Ob(\CCC)$. An object $B \in \Ob(\AAA)$ together with a
morphism $c : B \to C$ is a \emph{coreflection of $C$ in $\AAA$} if for every $A \in \Ob(\AAA)$ and every morphism
$f \in \hom_\CCC(A, C)$ there is a unique morphism $g \in \hom_\AAA(A, B)$ such that $c \cdot g = f$.

\begin{lemma}\label{rpemklei.lem.refl} \cite{masul-RPEMKLEI-MONAD}
  Let $\CCC$ be a locally small category such that all the morphisms in $\CCC$ are mono. Let $\AAA$ be a full
  subcategory of $\CCC$, let $A, B, D \in \Ob(\AAA)$ and $C \in \Ob(\CCC)$. If
  $C \longrightarrow (B)^A_{k,t}$ for some $k, t \ge 2$ and if
  $c : D \to C$ is a coreflection of $C$ in $\AAA$ then $D \longrightarrow (B)^A_{k,t}$.\qed
\end{lemma}

A category $\CCC$ has the \emph{(finite) Ramsey property} if
for every integer $k \ge 2$ and all $A, B \in \Ob(\CCC)$ there is a
$C \in \Ob(\CCC)$ such that $C \longrightarrow (B)^{A}_k$.

For $A \in \Ob(\CCC)$ let $t_\CCC(A)$ denote the least positive integer $n$ such that
for all $k \ge 2$ and all $B \in \Ob(\CCC)$ there exists a $C \in \Ob(\CCC)$ such that
$C \longrightarrow (B)^{A}_{k, n}$, if such an integer exists.
Otherwise put $t_\CCC(A) = \infty$. The number $t_\CCC(A)$ is referred to as the \emph{small Ramsey degree}
of $A$ in $\CCC$. A category $\CCC$ has the \emph{finite small Ramsey degrees} if $t_\CCC(A) < \infty$
for all $A \in \Ob(\CCC)$. Clearly, a category $\CCC$ has the Ramsey property if and only if $t_\CCC(A) = 1$
for all $A \in \Ob(\CCC)$. In this parlance the Finite Ramsey Theorem takes the following form.

\begin{theorem} [The Finite Ramsey Theorem~\cite{ramsey}]
  The category $\ChEmbFin$ has the Ramsey property.\qed
\end{theorem}

The Ramsey property for ordered structures implies the existence of finite small Ramsey degrees for
the corresponding unordered structures. This was first observed for categories of structures
in~\cite{dasilvabarbosa}, and generalized to arbitrary categories in~\cite{masul-dual-kpt}.

Let us outline the main tool we employ to obtain results of this form.
Following \cite{KPT,dasilvabarbosa,masul-dual-kpt} we say that
an \emph{expansion} of a category $\CCC$ is a category $\CCC^*$ together with a forgetful functor $U : \CCC^* \to \CCC$.
We shall generally follow the convention that $A, B, C, \dots$ denote objects from $\CCC$
while $\calA, \calB, \calC, \dots$ denote objects from $\CCC^*$.
Since $U$ is injective on hom-sets we may safely assume that
$\hom_{\CCC^*}(\calA, \calB) \subseteq \hom_\CCC(A, B)$ where $A = U(\calA)$, $B = U(\calB)$.
In particular, $\id_\calA = \id_A$ for $A = U(\calA)$. Moreover, it is safe to drop
subscripts $\CCC$ and $\CCC^*$ in $\hom_\CCC(A, B)$ and $\hom_{\CCC^*}(\calA, \calB)$, so we shall
simply write $\hom(A, B)$ and $\hom(\calA, \calB)$, respectively.
Let
$
  U^{-1}(A) = \{\calA \in \Ob(\CCC^*) : U(\calA) = A \}
$. Note that this is not necessarily a set.

An expansion $U : \CCC^* \to \CCC$ is \emph{reasonable} (cf.~\cite{KPT,masul-dual-kpt}) if
for every $e \in \hom(A, B)$ and every $\calA \in U^{-1}(A)$ there is a $\calB \in U^{-1}(B)$ such that
$e \in \hom(\calA, \calB)$.
An expansion $U : \CCC^* \to \CCC$ has \emph{unique restrictions} \cite{masul-dual-kpt} if
for every $\calB \in \Ob(\CCC^*)$ and every $e \in \hom(A, U(\calB))$ there is a \emph{unique} $\calA \in U^{-1}(A)$
such that $e \in \hom(\calA, \calB)$.
We denote this unique $\calA$ by $\restr \calB e$ and refer to it as the \emph{restriction of $\calB$ along~$e$}.

The following result was first proved for categories of structures in \cite{dasilvabarbosa},
and for general categories in \cite{masul-dual-kpt}.

\begin{theorem}\label{sbrd.thm.small1} \cite{dasilvabarbosa,masul-dual-kpt}
  Let $\CCC$ and $\CCC^*$ be locally small categories such that all the morphisms in $\CCC$ and $\CCC^*$ are mono.
  Let $U : \CCC^* \to \CCC$ be a reasonable expansion with unique restrictions. Then
  $
    t_{\CCC}(A) \le \sum_{\calA \in U^{-1}(A)} t_{\CCC^*}(\calA)
  $
  for all  $A \in \Ob(\CCC)$. Consequently, if $U^{-1}(A)$ is finite and
  $t_{\CCC^*}(\calA) < \infty$ for all $\calA \in U^{-1}(A)$ then $t_\CCC(A) < \infty$.

  In particular, if $U : \CCC^* \to \CCC$ is a reasonable expansion with unique restrictions such that $\CCC^*$ has the
  Ramsey property and $U^{-1}(A)$ is finite for all $A \in \Ob(\CCC)$ then $\CCC$ has finite small Ramsey degrees.\qed
\end{theorem}

Let $\CCC$ be a locally small category.
For $A, S \in \Ob(\CCC)$ let $T_\CCC(A, S)$ denote the least positive integer $n$ such that
$S \longrightarrow (S)^{A}_{k, n}$ for all $k \ge 2$, if such an integer exists.
Otherwise put $T_\CCC(A, S) = \infty$. The number $T_\CCC(A, S)$ is referred to as the \emph{big Ramsey degree}
of $A$ in $S$.
In this parlance the Infinite Ramsey Theorem takes the following form.

\begin{theorem}[The Infinite Ramsey Theorem~\cite{ramsey}] \label{rpemklei.thm.IRT}
  In the category $\ChEmb$ we have that
  $T(A, \omega) = 1$ for every finite chain~$A$.\qed
\end{theorem}

The following result was first proved for categories of structures in \cite{dasilvabarbosa},
and for general categories in \cite{masul-big-v-small}.

\begin{theorem}\label{rpemklei---big.v.small} (cf. \cite{dasilvabarbosa,masul-big-v-small})
  Let $\CCC$ and $\CCC^*$ be locally small categories such that all the morphisms in $\CCC$ and $\CCC^*$ are mono.
  Let $U : \CCC^* \to \CCC$ be an expansion with unique restrictions.
  For $A \in \Ob(\CCC)$, $S^* \in \Ob(\CCC^*)$ and $S = U(S^*)$, if $U^{-1}(A)$ is finite then
  $
    T_{\CCC}(A, S) \le \sum_{A^* \in U^{-1}(A)} T_{\CCC^*}(A^*, S^*)
  $.\qed
\end{theorem}

Our major tool for transporting the Ramsey property from one context to another
is to establish an adjunction-like relationship between the corresponding categories.

\begin{theorem}\label{rpemklei.thm.adjoints-rp} \cite{masul-scow}
  Right adjoints preserve the Ramsey property while left adjoints preserve the dual Ramsey property.
  More precisely, let $\BB$ and $\CCC$ be locally small categories and let $F : \BB \rightleftarrows \CCC : H$ be an adjunction.

  $(a)$ If $\CCC$ has the Ramsey property then so does $\BB$

  $(b)$ If $\BB$ has the dual Ramsey property then so does $\CCC$.\qed
\end{theorem}

\begin{theorem} \cite{masul-RPEMKLEI-MONAD}
  Let $\CCC$ be a locally small category, $(E, \delta, \epsilon)$ a comonad on $\CCC$, and let
  $\KK = \KK(E, \delta, \epsilon)$ and $\EM = \EM(E, \delta, \epsilon)$
  be the Kleisli category and the Eilenberg-Moore category, respectively, for the comonad. If $\CCC$ has the
  Ramsey property then so do both $\KK$ and $\EM$.
\end{theorem}

However, more is true in case of the Eilenberg-Moore construction.
We are now going to show that the Ramsey property
carries over from $\CCC$ to a more general context of coalgebras for functors with comultiplication,
which are straightforward weakenings of comonads. The proof relies on the following weakening of the notion of adjunction.

\begin{definition}\label{opos.def.PA} \cite{masul-preadj}
  Let $\BB$ and $\CCC$ be locally small categories. A pair of maps
  $
    F : \Ob(\BB) \rightleftarrows \Ob(\CCC) : H
  $
  is a \emph{pre-adjunction between $\BB$ and $\CCC$} provided there is a family of maps
  $
    \Phi_{X,Y} : \hom_\CCC(F(X), Y) \to \hom_\BB(X, H(Y))
  $
  indexed by the pairs $(X, Y) \in \Ob(\BB) \times \Ob(\CCC)$ and satisfying the following:
  \begin{itemize}
  \item[(PA)]
  for every $C \in \Ob(\CCC)$, every $A, B \in \Ob(\BB)$,
  every $u \in \hom_\CCC(F(B), C)$ and every $f \in \hom_\BB(A, B)$ there is a $v \in \hom_\CCC(F(A), F(B))$
  satisfying $\Phi_{B, C}(u) \cdot f = \Phi_{A, C}(u \cdot v)$.
  \end{itemize}
\end{definition}
\noindent
(Note that in a pre-adjunction $F$ and $H$ are \emph{not} required to be functors, just maps from the class of objects of one of the two
categories into the class of objects of the other category; also $\Phi$ is not required to be a natural isomorphism, just a family of
maps between hom-sets satisfying the requirement above.)

The proof of the following result is a straightforward generalization of \cite[Theorem 3.2]{masul-preadj}, but we shall
include it for completeness.

\begin{theorem}\label{opos.thm.main}
  Let $\BB$ and $\CCC$ be locally small categories and let $F : \Ob(\BB) \rightleftarrows \Ob(\CCC) : H$ be a pre-adjunction.
  \begin{itemize}
  \item[$(a)$]
    Let $t, k \ge 2$ be integers, let $A, B \in \Ob(\BB)$ and $C \in \Ob(\CCC)$.
    If $C \longrightarrow (F(B))^{F(A)}_{k,t}$ in $\CCC$ then $H(C) \longrightarrow (B)^A_{k,t}$ in $\BB$.
  \item[$(b)$]
    $t_\BB(A) \le t_\CCC(F(A))$ for all $A \in \Ob(\BB)$.
  \item[$(c)$]
    If $\CCC$ has the Ramsey property then so does $\BB$.\qed
  \end{itemize}
\end{theorem}
\begin{proof}
  $(a)$
  Assume that $A, B \in \Ob(\DD)$ and $k, t \ge 2$ satisfy $C \longrightarrow (F(B))^{F(A)}_{k, t}$
  and let us show that $H(C) \longrightarrow (B)^{A}_{k, t}$.
  Take any $\chi : \hom_\BB(A, H(C)) \to k$ and define
  $\chi' : \hom_\CCC(F(A), C) \to k$ by $\chi'(u) = \chi(\Phi_{A, C}(u))$.
  Since $C \longrightarrow (F(B))^{F(A)}_{k, t}$ there is a $u \in \hom_\CCC(F(B), C)$ such that
  \begin{equation}\label{opos.eq.2}
    |\chi'(u \cdot \hom_\CCC(F(A), F(B)))| \le t.
  \end{equation}
  Let us show that
  \begin{equation}\label{opos.eq.3}
    \chi(\Phi_{B, C}(u) \cdot \hom_\BB(A, B)) \subseteq \chi'(u \cdot \hom_\CCC(F(A), F(B))).
  \end{equation}
  Take any $f \in \hom_\BB(A, B)$. Since $F : \Ob(\BB) \rightleftarrows \Ob(\CCC) : H$ is a pre-adjunction, there is a
  $v \in \hom_\CCC(F(A), F(B))$ such that $\Phi_{B, C}(u) \cdot f = \Phi_{A, C}(u \cdot v)$. Then
  $
    \chi(\Phi_{B, C}(u) \cdot f) = \chi(\Phi_{A, C}(u \cdot v)) = \chi'(u \cdot v)
  $,
  which proves \eqref{opos.eq.3}. Finally, \eqref{opos.eq.3} together with \eqref{opos.eq.2}
  implies
  $$
    |\chi(\Phi_{B, C}(u) \cdot \hom_\BB(A, B))| \le t.
  $$

  $(b)$ and $(c)$ are now immediate.  
\end{proof}

The proof of the dual version of the theorem below is given in~\cite{masul-RPEMKLEI-MONAD}. Just as an illustration we provide the
proof of the ``direct'' version here.

\begin{theorem}\label{rpemklei.thm.ramsey} (cf.~\cite[Theorem 3.15]{masul-RPEMKLEI-MONAD})
  Let $\CCC$ be a locally small category, $E : \CCC \to \CCC$ a functor and $\delta : E \to EE$ a comultiplication for $E$.
  If $\CCC$ has the Ramsey property then so does every full subcategory of $\EM^w(E, \delta)$
  which contains all the $E$-coalgebras of the form $(E(C), \delta_C)$, $C \in \Ob(\CCC)$.
\end{theorem}
\begin{proof}
  Let $\BB$ be a full subcategory of $\EM^w(E, \delta)$ such that
  all the $E$-coalgebras of the form $(E(C), \delta_C)$, $C \in \Ob(\CCC)$, are in~$\BB$.
  By Theorem~\ref{opos.thm.main} in order to show that $\BB$ has the Ramsey property
  it suffices to construct a pre-adjunction $F : \Ob(\BB) \rightleftarrows \Ob(\CCC) : H$.
  For $(B, \beta) \in \Ob(\BB)$ put $F(B, \beta) = B$,
  for $C \in \Ob(\CCC)$ put $H(C) = (E(C), \delta_C)$ and
  for $u \in \hom_\CCC(B, C)$ put $\Phi_{(B, \beta), C}(u) = E(u) \cdot \beta$.

  Let us first show that the definition of $\Phi$ is correct by showing that for each
  $u \in \hom_\CCC(B, C)$ we have that $\Phi_{(B, \beta), C}(u)$ is a coalgebraic homomorphism
  from $(B, \beta)$ to $H(C)$, that is:
  \begin{center}
    \begin{tikzcd}
      B \arrow[rr, "E(u) \cdot \beta"] \arrow[d, "\beta"'] & & E(C) \arrow[d, "\delta_C"] \\
      E(B) \arrow[rr, "E(E(u) \cdot \beta)"'] & & EE(C)
    \end{tikzcd}
  \end{center}
  \noindent
  The following two diagrams commute because $(B, \beta)$ is a weak $E$-coalgebra and because
  $\delta : E \to EE$ is a natural transformation, respectively:
  \begin{center}
    \begin{tabular}{c@{\qquad}c}
      \begin{tikzcd}
        B \arrow[r, "\beta"] \arrow[d, "\beta"'] & E(B) \arrow[d, "E(\beta)"] \\
        E(B) \arrow[r, "\delta_B"'] & EE(B)
      \end{tikzcd}
      &
      \begin{tikzcd}
        E(B) \arrow[r, "E(u)"] \arrow[d, "\delta_B"'] & E(C) \arrow[d, "\delta_C"] \\
        EE(B) \arrow[r, "EE(u)"'] & EE(C)
      \end{tikzcd}
    \end{tabular}
  \end{center}
  It now easily follows that $EE(u) \cdot E(\beta) \cdot \beta = EE(u) \cdot \delta_B \cdot \beta
  = \delta_c \cdot E(u) \cdot \beta$.

  To complete the proof we still have to show that the condition (PA) of Definition~\ref{opos.def.PA}
  is satisfied.
  Take any $C \in \Ob(\CCC)$ and $(A, \alpha), (B, \beta) \in \Ob(\BB)$, take
  arbitrary $u \in \hom_\CCC(B, C)$ and an
  arbitrary coalgebraic homomorphism $f \in \hom_\BB((A, \alpha), (B, \beta))$.
  Then, $f \in \hom_\CCC(A, B)$ and
  $
    \Phi_{(B, \beta), C}(u) \cdot f
    = E(u) \cdot \beta \cdot f
    = E(u) \cdot E(f) \cdot \alpha
    = E(u \cdot f) \cdot \alpha  = \Phi_{(A, \alpha), C}(u \cdot f)
  $,
  having in mind that $\beta \cdot f = E(f) \cdot \alpha$ because $f$ is a coalgebraic homomorphism.
  This completes the proof.
\end{proof}

As a warm-up for the discussion that follows
we conclude the section by reproving a simple and well-known Ramsey result about finite trees \cite{deuber-finite-trees}
(see also \cite[p.~28]{leeb-vorlesungen}). Let $(A, f)$ be a monounary algebra.
We shall say that $r \in A$ is a \emph{root} if $f(r) = r$. A \emph{rooted forest} is a monounary algebra $(A, f)$ such that
for every $a \in A$ there is an $n \ge 0$ and a root $r$ such that $f^n(a) = r$. A \emph{rooted tree} is a rooted forest with
a single root. Let $\RfEmb$ denote the category of rooted forests with embeddings, and let
$\RfEmbFin$ denote the full subcategory of $\RfEmb$ spanned by its finite members.

An \emph{ordered rooted forest (resp.\ tree)} is a rooted forest (resp.\ tree) together with a linear ordering of its vertices.
An embedding between two ordered rooted forests is every monotone embedding of the underlying rooted forests.
Let $\OrfEmb$ denote the category of ordered rooted forests with embeddings, and let
$\OrfEmbFin$ denote the full subcategory of $\OrfEmb$ spanned by its finite members.
We are now going to prove the Ramsey property for the class of all finite ordered rooted forests
by representing them as $E$-coalgebras for a modification of the comonad presented in Example~\ref{rpemklei.ex.list-comonad}.

\begin{figure}
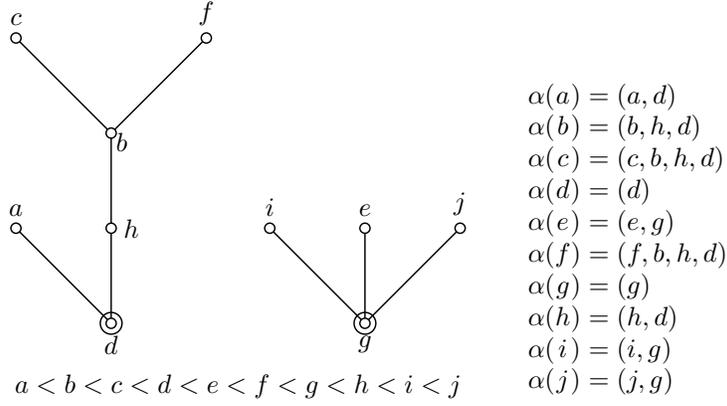

  \centering
\begin{pgfpicture}
  \pgfsetxvec{\pgfpoint{\acadpgfunit}{0pt}}
  \pgfsetyvec{\pgfpoint{0pt}{\acadpgfunit}}
  \pgfsetlinewidth{\acadpgflinewidth}
  \pgftransformshift{\pgfpointxy{25.0}{-75.0}}

  \begin{pgfscope}
    \pgfpathmoveto{\pgfpointxy{250.0}{250.0}}
    \pgfpathlineto{\pgfpointxy{100.0}{350.0}}
    \pgfusepath{stroke}
  \end{pgfscope}
  \begin{pgfscope}
    \pgfpathmoveto{\pgfpointxy{250.0}{250.0}}
    \pgfpathlineto{\pgfpointxy{250.0}{350.0}}
    \pgfusepath{stroke}
  \end{pgfscope}
  \begin{pgfscope}
    \pgfpathmoveto{\pgfpointxy{250.0}{350.0}}
    \pgfpathlineto{\pgfpointxy{250.0}{450.0}}
    \pgfusepath{stroke}
  \end{pgfscope}
  \begin{pgfscope}
    \pgfpathmoveto{\pgfpointxy{250.0}{450.0}}
    \pgfpathlineto{\pgfpointxy{100.0}{550.0}}
    \pgfusepath{stroke}
  \end{pgfscope}
  \begin{pgfscope}
    \pgfpathmoveto{\pgfpointxy{250.0}{450.0}}
    \pgfpathlineto{\pgfpointxy{400.0}{550.0}}
    \pgfusepath{stroke}
  \end{pgfscope}
  \begin{pgfscope}
    \pgfpathmoveto{\pgfpointxy{650.0}{250.0}}
    \pgfpathlineto{\pgfpointxy{500.0}{350.0}}
    \pgfusepath{stroke}
  \end{pgfscope}
  \begin{pgfscope}
    \pgfpathmoveto{\pgfpointxy{650.0}{250.0}}
    \pgfpathlineto{\pgfpointxy{650.0}{350.0}}
    \pgfusepath{stroke}
  \end{pgfscope}
  \begin{pgfscope}
    \pgfpathmoveto{\pgfpointxy{650.0}{250.0}}
    \pgfpathlineto{\pgfpointxy{800.0}{350.0}}
    \pgfusepath{stroke}
  \end{pgfscope}
  \begin{pgfscope}
    \pgfpathellipse{\pgfpointxy{250.0}{250.0}}{\pgfpointxy{17.1341}{0.0}}{\pgfpointxy{0.0}{17.1341}}
    \pgfusepath{stroke}
  \end{pgfscope}
  \begin{pgfscope}
    \pgfpathellipse{\pgfpointxy{650.0}{250.0}}{\pgfpointxy{17.1341}{0.0}}{\pgfpointxy{0.0}{17.1341}}
    \pgfusepath{stroke}
  \end{pgfscope}
  \begin{pgfscope}
    \pgfsetfillcolor{white}
    \pgfpathellipse{\pgfpointxy{250.0}{250.0}}{\pgfpointxy{8.0}{0.0}}{\pgfpointxy{0.0}{8.0}}
    \pgfusepath{fill,stroke}
  \end{pgfscope}
  \begin{pgfscope}
    \pgfsetfillcolor{white}
    \pgfpathellipse{\pgfpointxy{250.0}{350.0}}{\pgfpointxy{8.0}{0.0}}{\pgfpointxy{0.0}{8.0}}
    \pgfusepath{fill,stroke}
  \end{pgfscope}
  \begin{pgfscope}
    \pgfsetfillcolor{white}
    \pgfpathellipse{\pgfpointxy{100.0}{350.0}}{\pgfpointxy{8.0}{0.0}}{\pgfpointxy{0.0}{8.0}}
    \pgfusepath{fill,stroke}
  \end{pgfscope}
  \begin{pgfscope}
    \pgfsetfillcolor{white}
    \pgfpathellipse{\pgfpointxy{250.0}{450.0}}{\pgfpointxy{8.0}{0.0}}{\pgfpointxy{0.0}{8.0}}
    \pgfusepath{fill,stroke}
  \end{pgfscope}
  \begin{pgfscope}
    \pgfsetfillcolor{white}
    \pgfpathellipse{\pgfpointxy{100.0}{550.0}}{\pgfpointxy{8.0}{0.0}}{\pgfpointxy{0.0}{8.0}}
    \pgfusepath{fill,stroke}
  \end{pgfscope}
  \begin{pgfscope}
    \pgfsetfillcolor{white}
    \pgfpathellipse{\pgfpointxy{400.0}{550.0}}{\pgfpointxy{8.0}{0.0}}{\pgfpointxy{0.0}{8.0}}
    \pgfusepath{fill,stroke}
  \end{pgfscope}
  \begin{pgfscope}
    \pgfsetfillcolor{white}
    \pgfpathellipse{\pgfpointxy{650.0}{250.0}}{\pgfpointxy{8.0}{0.0}}{\pgfpointxy{0.0}{8.0}}
    \pgfusepath{fill,stroke}
  \end{pgfscope}
  \begin{pgfscope}
    \pgfsetfillcolor{white}
    \pgfpathellipse{\pgfpointxy{650.0}{350.0}}{\pgfpointxy{8.0}{0.0}}{\pgfpointxy{0.0}{8.0}}
    \pgfusepath{fill,stroke}
  \end{pgfscope}
  \begin{pgfscope}
    \pgfsetfillcolor{white}
    \pgfpathellipse{\pgfpointxy{500.0}{350.0}}{\pgfpointxy{8.0}{0.0}}{\pgfpointxy{0.0}{8.0}}
    \pgfusepath{fill,stroke}
  \end{pgfscope}
  \begin{pgfscope}
    \pgfsetfillcolor{white}
    \pgfpathellipse{\pgfpointxy{800.0}{350.0}}{\pgfpointxy{8.0}{0.0}}{\pgfpointxy{0.0}{8.0}}
    \pgfusepath{fill,stroke}
  \end{pgfscope}
  \pgftext[top,at={\pgfpointxy{250.0}{230.0}}]{$d$}
  \pgftext[bottom,at={\pgfpointxy{100.0}{370.0}}]{$a$}
  \pgftext[top,left,at={\pgfpointxy{258.0}{450.0}}]{$b$}
  \pgftext[bottom,at={\pgfpointxy{100.0}{570.0}}]{$c$}
  \pgftext[bottom,at={\pgfpointxy{650.0}{370.0}}]{$e$}
  \pgftext[bottom,at={\pgfpointxy{400.0}{570.0}}]{$f$}
  \pgftext[top,at={\pgfpointxy{650.0}{230.0}}]{$g$}
  \pgftext[left,at={\pgfpointxy{270.0}{350.0}}]{$h$}
  \pgftext[bottom,at={\pgfpointxy{500.0}{370.0}}]{$i$}
  \pgftext[bottom,at={\pgfpointxy{800.0}{370.0}}]{$j$}
  \pgftext[at={\pgfpointxy{450.0}{150.0}}]{$a<b<c<d<e<f<g<h<i<j$}
\end{pgfpicture}
  \qquad
  \begin{tabular}[b]{@{$\alpha($}c@{$) = \mathstrut$}l}
    $a$ & $(a,d)$ \\
    $b$ & $(b,h,d)$ \\
    $c$ & $(c,b,h,d)$ \\
    $d$ & $(d)$ \\
    $e$ & $(e,g)$ \\
    $f$ & $(f,b,h,d)$ \\
    $g$ & $(g)$ \\
    $h$ & $(h,d)$ \\
    $i$ & $(i,g)$ \\
    $j$ & $(j,g)$
  \end{tabular}
  \caption{An ordered rooted forest and it representation as an $E$-coalgebra}
  \label{rpemklei.fig.rootedforest}
\end{figure}

\begin{example}
  Fig.~\ref{rpemklei.fig.rootedforest} depicts an ordered rooted forest. Note that
  every ordered rooted forest can be represented by an $E$-coalgebra straightforwardly: if $(A, f)$ is an
  ordered rooted forest define $\alpha : A \to E(A)$ by $\alpha(a) = \mathstrut$the unique path from $a$ to the root of its connected
  component. 
\end{example}

\begin{theorem}\label{rpemklei.thm.orf-fin} (\cite{deuber-finite-trees}; cf.\ \cite[p.~28]{leeb-vorlesungen})
  The category $\OrfEmbFin$ of finite ordered rooted forests and embeddings has the Ramsey property.
\end{theorem}
\begin{proof}
  For a set $A$ let $A^\dagger$ denote the set of all finite sequences of pairwise distinct elements of $A$.
  (Note that if $A$ is finite then so is $A^\dagger$.) Let $E : \ChEmbFin \to \ChEmbFin$ be the functor such that
  $E(A, \Boxed{\lt A}) = (A^\dagger, \Boxed{\lex A})$ and $E(f)(a_1, \dots, a_k) = (f(a_1), \dots, f(a_k))$.
  If $f : (A, \Boxed{\lt A}) \to (B, \Boxed{\lt B})$ is an embedding then $E(f) : (A^\dagger, \Boxed{\lex A}) \to (B^\dagger, \Boxed{\lex B})$
  is also an embedding, so the definition of $E$ is correct. Finally, let us define comultiplication $\delta : E \to EE$ as
  in Example~\ref{rpemklei.ex.list-comonad}:
  $$
     \delta_A(a_1, \dots, a_n) = \big(
       (a_1, a_2, \dots, a_n), (a_2, \dots, a_n), \dots, (a_n)
     \big)
  $$
  It is obvious that, given a chain $(A, \Boxed{\lt A})$, the map $\delta_A : E(A, \Boxed{\lt A}) \to EE(A, \Boxed{\lt A})$ is an embedding
  of chains, so the definition of $\delta$ is also correct.

  It follows immediately from the construction that each $\alpha$ is an
  embedding between the corresponding chains, and hence, a morphism in $\ChEmbFin$. Moreover, every coalgebra constructed
  in this manner is a weak Eilenberg-Moore coalgebra for $E$ and $\delta$. It is also easy to see that every embedding between
  two ordered rooted forests is clearly a coalgebraic homomorphism and vice versa.
  Moreover, every cofree Eilenberg-Moore $E$-coalgebra is a $E$-coalgebra representation of an ordered rooted forest.
  Therefore, $\OrfEmbFin$ is a full subcategory of $\EE$, the category whose objects are weak Eilenberg-Moore $E$-coalgebras
  and morphisms are coalgebraic homomorphisms, and contains all the cofree Eilenberg-Moore $E$-coalgebras.
  The claim now follows by Theorem~\ref{rpemklei.thm.ramsey}.
\end{proof}

\begin{corol}\label{rpemklei.cor.rooted-forests}
  The category $\RfEmbFin$ of finite rooted forests and embeddings has finite Ramsey degrees.
\end{corol}
\begin{proof}
  Let $U : \OrfEmbFin \to \RfEmbFin$ be the forgetful functor that forgets the order.
  It is now easy to see that $U$ is a reasonable expansion with unique restrictions. Since
  $\OrfEmbFin$ has the Ramsey property and $U^{-1}(A)$ is finite for every finite rooted forest $A$
  (because there are only finitely many linear orders on a finite set) Theorem~\ref{sbrd.thm.small1} implies that
  $\RfEmbFin$ has finite Ramsey degrees.
\end{proof}

\section{Ramsey properties of $M$-sets}
\label{rpemklei.sec.unary-algs}

In this section we apply the proof strategies outlined in Section~\ref{rpemklei.sec.RP} by proving
several Ramsey results for categories of $M$-sets. In 2016 Soki\'c proved that
for a finite monoid $M$ the class of all ordered finite $M$-sets has the Ramsey property~\cite{sokic-unary-functions}.
Using a completely different strategy, in this section we prove that for any monoid $M$ (finite or infinite) the category
of all finite ordered $M$-sets with embeddings has the Ramsey property, and from that conclude that
(unordered) finite $M$-sets have finite small Ramsey degrees. Moreover, we prove that for any monoid $M$ (finite or infinite)
finite ordered $M$-sets have finite big Ramsey degrees in the ordered cofree $M$-set $\hat \calE(\omega)$ on $\omega$ generators
and again infer the corresponding result for the unordered case.

Our proof mimics the proof of the fact that the category of weak Eilenberg-Moore coalgebras for a comonad
has the Ramsey property (Theorem~\ref{rpemklei.thm.ramsey}). Unfortunately, we are unable to apply Theorem~\ref{rpemklei.thm.ramsey} directly
because the construction in this section relies on the comonad $E(A) = A^M$ of Example~\ref{rpemklei.ex.mon-act} which has the unpleasant property
that $E(A)$ is infinite in case $M$ is an infinite monoid; and treating the case of infinite monoids is the key motivation for the paper.

We shall bypass this issue using the following compactness argument which was first proved for
categories of structures in~\cite{mu-pon} (see also~\cite{vanthe-more}), and for general categories in \cite{masul-dual-kpt}.
An $F \in \Ob(\CCC)$ is \emph{weakly locally finite for $\BB$} (cf.\ locally finite in~\cite{masul-dual-kpt}) if
for every $A, B \in \Ob(\BB)$ and every $e \in \hom(A, F)$, $f \in \hom(B, F)$ there exist $D \in \Ob(\BB)$,
$r \in \hom(D, F)$, $p \in \hom(A, D)$ and $q \in \hom(B, D)$ such that $r \cdot p = e$ and $r \cdot q = f$.

\begin{lemma}\label{akpt.lem.ramseyF} \cite{mu-pon,vanthe-more,masul-dual-kpt}
  Let $\CCC$ be a locally small category whose morphisms are mono, let $\BB$ be a full subcategory of $\CCC$
  such that $\hom_\BB(A, B)$ is finite for all $A, B \in \Ob(\BB)$, and let $F \in \Ob(\CCC)$ be universal and weakly locally finite
  for $\BB$. Then for all $A \in \Ob(\BB)$ and $t \ge 2$:
    $t_{\BB}(A) \le t$ if and only if
    $F \longrightarrow (B)^A_{k, t}$ for all $k \ge 2$ and all $B \in \Ob(\BB)$ such that $A \overset\BB\longrightarrow B$.\qed
\end{lemma}

Let $M$ be a monoid (finite or infinite). As we have already seen in Example~\ref{rpemklei.ex.mon-act}
every $M$-set can be represented by an Eilenberg-Moore coalgebra for the comonad
$E : \Set \to \Set$ defined by $E(A) = A^M$ on objects and by $E(f) : A^M \to B^M : h \mapsto f \circ h$
on morphisms so that $\SetEmb(M) \cong \EM_\emb(E, \delta, \epsilon)$.
In the proofs below we shall move freely between the two representations of $M$-sets.

Let us now upgrade $E$ and $\delta$ to $\ChEmb$. Take an arbitrary but fixed well-ordering of $M$ such that
$\min M = 1$. Recall that for every chain $(X, \Boxed<)$ the set $X^{M}$ can be ordered lexicographically as follows:
for $f, g \in X^{M}$ such that $f \ne g$ let
$f \lex{X} g \text{ iff } f(v) < g(v)$, where $v = \min \{w \in M : f(w) \ne g(w)\}$.
For a chain $(X, \Boxed<)$ let
$
  \hat E(X, \Boxed<) = (X^{M}, \Boxed{\lex X})
$.
This is how $\hat E$ acts on objects.
For an embedding $h : (X, \Boxed<) \to (Y, \Boxed<)$ define $\hat E(h) : \hat E(X, \Boxed<) \to \hat E(Y, \Boxed<)$ as
$
  \hat E(h)(f) = E(h)(f) = h \cdot f
$.
To see that $\hat E$ is well-defined take $f, g \in X^{M}$ such that $f \lex{X} g$.
Let $v = \min \{w \in M : f(w) \ne g(w)\}$. Then $f(w) < g(w)$ for $w < v$ and $f(v) = g(v)$. Since
$h$ is an embedding we immediately get that $h(f(w)) < h(g(w))$ for $w < v$ and $h(f(v)) = h(g(v))$, whence
$\hat E(h)(f) \lex{X} \hat E(h)(g)$. In other words, $\hat E(h)$ is an embedding $(X^{M}, \Boxed{\lex X}) \to (Y^{M}, \Boxed{\lex Y})$.

Following Example~\ref{rpemklei.ex.mon-act} let us define comultiplication $\hat\delta_{(X, \Boxed<)} : \hat E(X, \Boxed<) \to \hat E \hat E(X, \Boxed<)$
by
$
  \hat \delta_{(X, \Boxed<)}(h)(v)(w) = h(wv)
$.
Let us show that the definition is correct, that is, that $\hat \delta_{(X, \Boxed<)}$ is an embedding.
Note that $\hat \delta_{(X, \Boxed<)}(h)(1) = h$. Take any $f, g \in \hat E(X, \Boxed<)$ such that $f \lex{X} g$. Then
$\hat \delta_{(X, \Boxed<)}(f)(1) < \hat \delta_{(X, \Boxed<)}(g)(1)$, whence $\hat \delta_{(X, \Boxed<)}(f) \lex{\hat E(X, \Boxed<)} \hat \delta_{(X, \Boxed<)}(g)$
because the ordering of $\hat E \hat E(X, \Boxed<) = (\hat E(X, \Boxed<))^{M}$ is lexicographic and $M$ is well-ordered
so that $1 = \min M$.

\begin{lemma}\label{rpemklei.lem.representation-of-M-sets}
  Every ordered $M$-set $\calA = (A, \alpha', \Boxed<)$ where $<$ is a linear
  order on $A$, can be represented by a weak Eilenberg-Moore $\hat E$-coalgebra $((A,\Boxed<), \alpha)$, where the
  structure map $\alpha : (A, \Boxed<) \to \hat E(A, \Boxed<)$ is defined by $\alpha(a)(g) = \alpha'(g, a)$.
\end{lemma}
\begin{proof}
  Clearly, we only have to check that $\alpha$ is an embedding. Take $a_1, a_2 \in A$ such that $a_1 < a_2$.
  Then $\alpha(a_1)(1) = a_1$ and $\alpha(a_2)(1) = a_2$, whence $\alpha(a_1) \lex{A} \alpha(a_2)$ because the ordering of
  $\hat E(A, \Boxed<) = (A, \Boxed<)^{M}$ is lexicographic and $1 = \min M$.
\end{proof}

Let $\OsetEmb(M)$ denote the category of ordered $M$-sets and embeddings understood
as a full subcategory of $\EM^w_\emb(\hat E, \hat \delta)$ using the correspondence given in the above lemma.

\begin{lemma}\label{rpemklei.lem.Unary-1}
  Taking $\EM^w_\emb(\hat E, \hat \delta)$ as the ambient category,
  for all $\calU, \calV \in \OsetEmbFin(M)$ and all $k \ge 2$ we have that
  $(\hat E(\omega), \hat\delta_{\omega}) \longrightarrow (\calV)^\calU_k$.
\end{lemma}
\begin{proof}
  Let $\BB = \OsetEmb(M)$ and let $F : \Ob(\BB) \rightleftarrows \Ob(\ChEmb) : H$ be a pre-adjunction
  constructed as in the proof of Theorem~\ref{rpemklei.thm.ramsey}:
  for $\calB = ((B, \Boxed{<}), \beta) \in \Ob(\BB)$ put $F(\calB) = (B, \Boxed{<})$,
  for $(C, \Boxed{<}) \in \Ob(\ChEmb)$ put $H(C, \Boxed{<}) = (\hat E(C, \Boxed{<}), \hat\delta_{(C, \Boxed{<})})$ and
  for $u \in \hom_{\ChEmb}((B, \Boxed{<}), (C, \Boxed{<}))$ put $\Phi_{\calB, (C, \Boxed{<})}(u) = \hat E(u) \cdot \beta$.

  Take any $\calU, \calV \in \OsetEmbFin(M)$ and any $k \ge 2$. Since $\ChEmbFin$ has the Ramsey property
  there is a finite chain $(W, \Boxed{<})$ such that $(W, \Boxed{<}) \longrightarrow (F(\calV))^{F(\calU)}_k$.
  By Theorem~\ref{opos.thm.main}~$(a)$ it then follows that $H(W, \Boxed{<}) = (\hat E(W, \Boxed{<}), \hat\delta_{(W, \Boxed{<})})
  \longrightarrow (\calV)^{\calU}_k$.
  Now, take any embedding $f : (W, \Boxed{<}) \to \omega$. The fact that $\hat \delta$ is natural
  yields that $\hat E(f) : \hat E(W, \Boxed{<}) \to \hat E(\omega)$
  is a morphism in $\BB$ from $(\hat E(W, \Boxed{<}), \hat\delta_{(W, \Boxed{<})})$ to $(\hat E(\omega), \delta_{\omega})$.
  Lemma~\ref{rpemklei.lem.C-D} now
  ensures that $(\hat E(\omega), \hat\delta_\omega) \longrightarrow (\calV)^{\calU}_k$.
\end{proof}

\begin{lemma}\label{rpemklei.lem.Unary-univ-loc-fin}
  Taking $\EM^w_\emb(\hat E, \hat\delta)$ as the ambient category,
  $(\hat E(\omega), \hat\delta_\omega)$ is universal and weakly locally finite for $\OsetEmbFin(M)$.
\end{lemma}
\begin{proof}
  Let $\BB = \OsetEmbFin(M)$.
  To see that $(\hat E(\omega), \hat\delta_\omega)$ is universal for $\BB$ take any $\calB = ((B, \Boxed{<}), \beta) \in \Ob(\BB)$
  and any embedding $f : (B, \Boxed{<}) \to (\omega, \Boxed<)$.
  Then the square on the right commutes because $\hat\delta$ is natural,
  while the square on the left commutes because $\calB$ is a weak Eilenberg-Moore $\hat E$-coalgebra:
  \begin{center}
    \begin{tikzcd}
          (B, \Boxed{<})         \arrow[rr, "\beta"]      \arrow[d, "\beta"']
      & & \hat E(B, \Boxed{<})   \arrow[rr, "\hat E(f)"]  \arrow[d, "\hat\delta_{(B,\Boxed<)}"]
      & & \hat E(\omega)                                  \arrow[d, "\hat\delta_{\omega}"]
    \\
          \hat E(B, \Boxed{<})          \arrow[rr, "\hat E(\beta)"']
      & & \hat E\hat E(B, \Boxed{<})    \arrow[rr, "\hat E\hat E(f)"']
      & & \hat E\hat E(\omega)
    \end{tikzcd}
  \end{center}
  Therefore, $\hat E(f) \cdot \beta \in \hom(\calB, (\hat E(\omega), \hat\delta_\omega))$.

  To see that $(\hat E(\omega), \hat\delta_\omega)$ is weakly locally finite for $\BB$ take any
  $\calA = ((A, \Boxed{<}), \alpha)$ and $\calB = ((B, \Boxed{<}), \beta)$ in $\Ob(\BB)$ and
  arbitrary morphisms $f : \calA \to (\hat E(\omega), \hat\delta_\omega)$ and
  $g : \calB \to (\hat E(\omega), \hat\delta_\omega)$. Then
  $f : (A, \alpha) \to (E(\omega), \delta_\omega)$ and $g : (B, \beta) \to (E(\omega), \delta_\omega)$
  are embeddings in $\SetEmb(M)$ of (unordered) $M$-sets and embeddings.
  So, $f(A)$ and $g(B)$ are carriers of two finite subcoalgebras of $(E(\omega), \delta_\omega)$.
  It is easy to see that $C = f(A) \cup g(B)$ is then also a carrier of a finite subcoalgebra of $(E(\omega), \delta_\omega)$,
  so let $\gamma : C \to E(C)$ be the structure map that turns $C$ into a subcoalgebra $(C, \gamma)$ of $(E(\omega), \delta_\omega)$.
  Therefore, the following diagram commutes in $\SetEmb(M)$, where $f_C : (A, \alpha) \to (C, \gamma)$ and
  $g_C : (B, \beta) \to (C, \gamma)$ are codomain restrictions of $f$ and $g$, respectively:
  \begin{center}
    \begin{tikzcd}
      (E(\omega), \delta_\omega) \arrow[rr, leftarrow, "\supset"] \arrow[d, leftarrow, "f"'] \arrow[drr, leftarrow, "g"', near end] & & (C, \gamma) \arrow[dll, leftarrow, "f_C", near end] \arrow[d, leftarrow, "g_C"] \\
      (A, \alpha) & & (B, \beta)
    \end{tikzcd}
  \end{center}
  Finally, let us order $C$ by restricting the linear ordering of $\hat E(\omega)$ to $C$. Then $((C, \Boxed<), \gamma) \in \Ob(\BB)$
  and all the morphisms in the diagram above are embeddings. This concludes the proof that $(\hat E(\omega), \hat\delta_\omega)$
  is locally finite for $\BB$.
\end{proof}

\begin{theorem}\label{rpemklei.thm.un-ord-rp}
  $(a)$ Let $M$ be an arbitrary monoid (finite or infinite).
  Then the category $\OsetEmbFin(M)$ of ordered finite $M$-sets and embeddings
  has the Ramsey property.

  $(b)$ For every group $G$ the category $\OsetEmbFin(G)$ of ordered finite $G$-sets and embeddings has the Ramsey property
        (for finite groups this was proved in \cite{sokic-unary-functions}).

  $(c)$ For every unary algebraic language $\Omega$ the category $\OalgEmbFin(\Omega)$ of ordered finite $\Omega$-algebras and embeddings has the
        Ramsey property (for finite unary languages $\Omega$ this was proved in \cite{sokic-unary-functions}).
\end{theorem}
\begin{proof}
  $(a)$ Fix a well-ordering of $M$ such that $1 = \min M$.
  Using this well-ordering construct the functor $\hat E : \ChEmb \to \ChEmb$ and the comultiplication
  $\hat \delta : \hat E \to \hat E \hat E$ as above. For notational convenience let $\CCC = \OsetEmbFin(M)$.
  Recall that $\CCC$ is a category of finite weak Eilenberg-Moore $\hat E$-coalgebras (see Lemma~\ref{rpemklei.lem.representation-of-M-sets}
  and the remark that follows). Take $\EM = \EM^w_\emb(\hat E, \hat \delta)$ as the ambient category.
  We have seen in Lemma~\ref{rpemklei.lem.Unary-univ-loc-fin} that
  $(\hat E(\omega), \hat\delta_{\omega})$ is universal and locally finite for $\CCC$.
  Lemma~\ref{rpemklei.lem.Unary-1} shows that
  for all $\calA, \calB \in \Ob(\CCC)$ and all $k \ge 2$ we have that $(\hat E(\omega), \hat\delta_{\omega}) \longrightarrow (\calB)^\calA_k$.
  Therefore, by Lemma~\ref{akpt.lem.ramseyF} we have that $t_{\CCC}(\calA) = 1$ for all $\calA \in \Ob(\CCC)$.
  This is just another way of saying that $\CCC$ has the Ramsey property.
  
  $(b)$ is a special case of $(a)$ where $M = G$ is a group.
  
  $(c)$ is a special case of $(a)$ where $M = \Omega^*$.
\end{proof}

\begin{corol}\label{rpemklei.cor.fsrd-G-set}
  For every monoid $M$ the category $\SetEmbFin(M)$ of finite $M$-sets and embeddings
  has finite small Ramsey degrees. In particular, 
  \begin{itemize}
  \item for every group $G$ the category $\SetEmbFin(G)$ of finite $G$-sets and embeddings has small Ramsey degrees; and
  \item for every unary algebraic language $\Omega$ the category $\AlgEmbFin(\Omega)$ of finite $\Omega$-algebras and embeddings has small Ramsey degrees.
  \end{itemize}
\end{corol}
\begin{proof}
  Let $U : \OsetEmbFin(M) \to \SetEmbFin(M)$ be the forgetful functor that forgets the order.
  It is now easy to see that $U$ is a reasonable expansion with unique restrictions. Since
  $\OsetEmbFin(M)$ has the Ramsey property and $U^{-1}(\calA)$ is finite for every finite $M$-set $\calA$
  (because there are only finitely many linear orders on a finite set) Theorem~\ref{sbrd.thm.small1} implies that
  $\SetEmbFin(M)$ has finite small Ramsey degrees.
\end{proof}

\begin{theorem}\label{rpemklei.thm.main-G-set-order}
  Let $M$ be an arbitrary monoid (finite or infinite).
  There exists an ordering $\hat\calE(\omega)$ of the cofree $M$-set $\calE(\omega)$ on $\omega$ generators
  such that every finite ordered $M$-set has finite big Ramsey degree in $\hat\calE(\omega)$.
  More precisely, for every finite ordered $M$-set $\calA = (A, \alpha, \Boxed<)$ we have
  $T(\calA, \hat\calE(\omega)) \le 2^{|A|-1}$ in $\OsetEmb(M)$.
\end{theorem}
\begin{proof}
  Fix a well-ordering of $M$ such that $1 = \min M$.
  Using this well-ordering construct the functor $\hat E : \ChEmb \to \ChEmb$ and the comultiplication
  $\hat \delta : \hat E \to \hat E \hat E$ as above. For notational convenience let $\CCC = \OsetEmbFin(M)$.
  Recall that $\CCC$ is a category of finite weak Eilenberg-Moore $\hat E$-coalgebras (see Lemma~\ref{rpemklei.lem.representation-of-M-sets}
  and the remark that follows). Let $\hat\calE(\omega) = (\hat E(\omega), \hat\delta_{\omega})$.

  Take any $\calA = ((A, \Boxed<), \alpha) \in \CCC^\fin$ where $(A, \Boxed<) = \{a_1 < a_2 < \dots < a_s\}$, $s = |A|$,
  and let us show that $T_\CCC(\calA, \hat\calE(\omega)) \le 2^{s-1}$.
  Let $\chi : \hom_\CCC(\calA, \hat\calE(\omega)) \to k$, $k \ge 2$, be an arbitrary coloring.

  Let $n = 2^{s-1}$ and let $A_1, \dots, A_n$ be all the subchains of $A$ that contain $a_1$.
  As a notational convenience, let
  $
    R = \hom_\CCC(\calA, \hat\calE(\omega))
  $
  and let
  $S_i = \hom_{\ChEmb}(A_i, \omega)$, $1 \le i \le n$.
  Take any $f \in R$ and let $f(a_i) = h_i \in \omega^M$, $1 \le i \le s$.
  Since $f$ is an embedding and $\omega^M$ is ordered lexicographically, we have that $h_1(1) \le h_2(1) \le \dots \le h_s(1)$.
  Define an equivalence relation $\rho$ on $A$ by $(a_i, a_j) \in \rho$ iff $h_i(1) = h_j(1)$ and let
  $
    A/\rho = \{B_1, \dots, B_m\}
  $
  where $\min B_1 < \dots < \min B_m$. Let $\min B_i = a_{j_i}$, $1 \le i \le m$ and note that $j_1 = 1$. Therefore,
  $\{a_{j_1}, \dots, a_{j_m}\}$ is a subset of $A$ that contains $a_1$, say, $A_\ell = \{a_{j_1}, \dots, a_{j_m}\}$.
  Finally, $f^* : A_\ell \to \omega : a_{j_i} \mapsto h_{j_i}(1)$ is clearly an embedding,
  see Fig.~\ref{rpemklei.fig.snopovi}.
  This defines a mapping
  $
    \pi : R \to \bigcup_{\ell=1}^n S_\ell : f \mapsto f^*
  $.

\begin{figure}
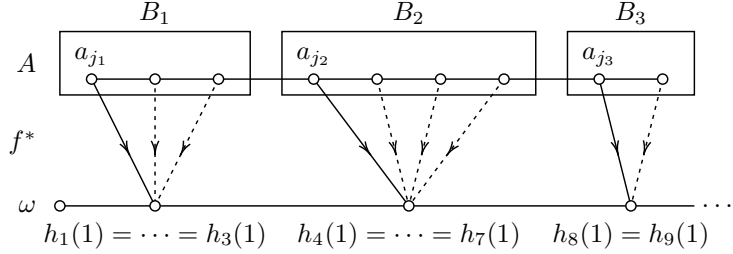

  \centering
\begin{pgfpicture}
  \pgfsetxvec{\pgfpoint{\acadpgfunit}{0pt}}
  \pgfsetyvec{\pgfpoint{0pt}{\acadpgfunit}}
  \pgfsetlinewidth{\acadpgflinewidth}
  \pgftransformshift{\pgfpointxy{100.0}{-50.0}}

  \begin{pgfscope}
    \pgfpathmoveto{\pgfpointxy{100.0}{200.0}}
    \pgfpathlineto{\pgfpointxy{1100.0}{200.0}}
    \pgfusepath{stroke}
  \end{pgfscope}
  \begin{pgfscope}
    \pgfpathmoveto{\pgfpointxy{400.0}{375.0}}
    \pgfpathlineto{\pgfpointxy{100.0}{375.0}}
    \pgfpathlineto{\pgfpointxy{100.0}{475.0}}
    \pgfpathlineto{\pgfpointxy{400.0}{475.0}}
    \pgfpathclose
    \pgfusepath{stroke}
  \end{pgfscope}
  \begin{pgfscope}
    \pgfpathmoveto{\pgfpointxy{850.0}{375.0}}
    \pgfpathlineto{\pgfpointxy{450.0}{375.0}}
    \pgfpathlineto{\pgfpointxy{450.0}{475.0}}
    \pgfpathlineto{\pgfpointxy{850.0}{475.0}}
    \pgfpathclose
    \pgfusepath{stroke}
  \end{pgfscope}
  \begin{pgfscope}
    \pgfpathmoveto{\pgfpointxy{1100.0}{375.0}}
    \pgfpathlineto{\pgfpointxy{900.0}{375.0}}
    \pgfpathlineto{\pgfpointxy{900.0}{475.0}}
    \pgfpathlineto{\pgfpointxy{1100.0}{475.0}}
    \pgfpathclose
    \pgfusepath{stroke}
  \end{pgfscope}
  \begin{pgfscope}
    \pgfpathmoveto{\pgfpointxy{150.0}{400.0}}
    \pgfpathlineto{\pgfpointxy{1050.0}{400.0}}
    \pgfusepath{stroke}
  \end{pgfscope}
  \begin{pgfscope}
    \pgfpathmoveto{\pgfpointxy{150.0}{400.0}}
    \pgfpathlineto{\pgfpointxy{250.0}{200.0}}
    \pgfusepath{stroke}
  \end{pgfscope}
  \begin{pgfscope}
    \pgfsetdash{{1.5pt}{2pt}}{0pt}
    \pgfpathmoveto{\pgfpointxy{250.0}{400.0}}
    \pgfpathlineto{\pgfpointxy{250.0}{200.0}}
    \pgfusepath{stroke}
  \end{pgfscope}
  \begin{pgfscope}
    \pgfsetdash{{1.5pt}{2pt}}{0pt}
    \pgfpathmoveto{\pgfpointxy{350.0}{400.0}}
    \pgfpathlineto{\pgfpointxy{250.0}{200.0}}
    \pgfusepath{stroke}
  \end{pgfscope}
  \begin{pgfscope}
    \pgfpathmoveto{\pgfpointxy{500.0}{400.0}}
    \pgfpathlineto{\pgfpointxy{650.0}{200.0}}
    \pgfusepath{stroke}
  \end{pgfscope}
  \begin{pgfscope}
    \pgfsetdash{{1.5pt}{2pt}}{0pt}
    \pgfpathmoveto{\pgfpointxy{600.0}{400.0}}
    \pgfpathlineto{\pgfpointxy{650.0}{200.0}}
    \pgfusepath{stroke}
  \end{pgfscope}
  \begin{pgfscope}
    \pgfsetdash{{1.5pt}{2pt}}{0pt}
    \pgfpathmoveto{\pgfpointxy{700.0}{400.0}}
    \pgfpathlineto{\pgfpointxy{650.0}{200.0}}
    \pgfusepath{stroke}
  \end{pgfscope}
  \begin{pgfscope}
    \pgfsetdash{{1.5pt}{2pt}}{0pt}
    \pgfpathmoveto{\pgfpointxy{800.0}{400.0}}
    \pgfpathlineto{\pgfpointxy{650.0}{200.0}}
    \pgfusepath{stroke}
  \end{pgfscope}
  \begin{pgfscope}
    \pgfpathmoveto{\pgfpointxy{950.0}{400.0}}
    \pgfpathlineto{\pgfpointxy{1000.0}{200.0}}
    \pgfusepath{stroke}
  \end{pgfscope}
  \begin{pgfscope}
    \pgfsetdash{{1.5pt}{2pt}}{0pt}
    \pgfpathmoveto{\pgfpointxy{1050.0}{400.0}}
    \pgfpathlineto{\pgfpointxy{1000.0}{200.0}}
    \pgfusepath{stroke}
  \end{pgfscope}
  \begin{pgfscope}
    \pgfpathmoveto{\pgfpointxy{1031.31}{300.366}}
    \pgfpatharcaxes{135.964}{165.964}{\pgfpointxy{45.0}{0.0}}{\pgfpointxy{0.0}{45.0}}
    \pgfusepath{stroke}
  \end{pgfscope}
  \begin{pgfscope}
    \pgfpathmoveto{\pgfpointxy{1020.0}{280.0}}
    \pgfpatharcaxes{-14.0362}{15.9638}{\pgfpointxy{45.0}{0.0}}{\pgfpointxy{0.0}{45.0}}
    \pgfusepath{stroke}
  \end{pgfscope}
  \begin{pgfscope}
    \pgfpathmoveto{\pgfpointxy{980.392}{303.29}}
    \pgfpatharcaxes{164.036}{194.036}{\pgfpointxy{45.0}{0.0}}{\pgfpointxy{0.0}{45.0}}
    \pgfusepath{stroke}
  \end{pgfscope}
  \begin{pgfscope}
    \pgfpathmoveto{\pgfpointxy{980.0}{280.0}}
    \pgfpatharcaxes{14.0362}{44.0362}{\pgfpointxy{45.0}{0.0}}{\pgfpointxy{0.0}{45.0}}
    \pgfusepath{stroke}
  \end{pgfscope}
  \begin{pgfscope}
    \pgfpathmoveto{\pgfpointxy{728.323}{294.383}}
    \pgfpatharcaxes{113.13}{143.13}{\pgfpointxy{45.0}{0.0}}{\pgfpointxy{0.0}{45.0}}
    \pgfusepath{stroke}
  \end{pgfscope}
  \begin{pgfscope}
    \pgfpathmoveto{\pgfpointxy{710.0}{280.0}}
    \pgfpatharcaxes{323.13}{353.13}{\pgfpointxy{45.0}{0.0}}{\pgfpointxy{0.0}{45.0}}
    \pgfusepath{stroke}
  \end{pgfscope}
  \begin{pgfscope}
    \pgfpathmoveto{\pgfpointxy{681.306}{300.366}}
    \pgfpatharcaxes{135.964}{165.964}{\pgfpointxy{45.0}{0.0}}{\pgfpointxy{0.0}{45.0}}
    \pgfusepath{stroke}
  \end{pgfscope}
  \begin{pgfscope}
    \pgfpathmoveto{\pgfpointxy{670.0}{280.0}}
    \pgfpatharcaxes{-14.0362}{15.9638}{\pgfpointxy{45.0}{0.0}}{\pgfpointxy{0.0}{45.0}}
    \pgfusepath{stroke}
  \end{pgfscope}
  \begin{pgfscope}
    \pgfpathmoveto{\pgfpointxy{630.392}{303.29}}
    \pgfpatharcaxes{164.036}{194.036}{\pgfpointxy{45.0}{0.0}}{\pgfpointxy{0.0}{45.0}}
    \pgfusepath{stroke}
  \end{pgfscope}
  \begin{pgfscope}
    \pgfpathmoveto{\pgfpointxy{630.0}{280.0}}
    \pgfpatharcaxes{14.0362}{44.0362}{\pgfpointxy{45.0}{0.0}}{\pgfpointxy{0.0}{45.0}}
    \pgfusepath{stroke}
  \end{pgfscope}
  \begin{pgfscope}
    \pgfpathmoveto{\pgfpointxy{581.323}{301.617}}
    \pgfpatharcaxes{186.87}{216.87}{\pgfpointxy{45.0}{0.0}}{\pgfpointxy{0.0}{45.0}}
    \pgfusepath{stroke}
  \end{pgfscope}
  \begin{pgfscope}
    \pgfpathmoveto{\pgfpointxy{590.0}{280.0}}
    \pgfpatharcaxes{36.8699}{66.8699}{\pgfpointxy{45.0}{0.0}}{\pgfpointxy{0.0}{45.0}}
    \pgfusepath{stroke}
  \end{pgfscope}
  \begin{pgfscope}
    \pgfpathmoveto{\pgfpointxy{305.455}{297.428}}
    \pgfpatharcaxes{123.435}{153.435}{\pgfpointxy{45.0}{0.0}}{\pgfpointxy{0.0}{45.0}}
    \pgfusepath{stroke}
  \end{pgfscope}
  \begin{pgfscope}
    \pgfpathmoveto{\pgfpointxy{290.0}{280.0}}
    \pgfpatharcaxes{-26.5651}{3.43495}{\pgfpointxy{45.0}{0.0}}{\pgfpointxy{0.0}{45.0}}
    \pgfusepath{stroke}
  \end{pgfscope}
  \begin{pgfscope}
    \pgfpathmoveto{\pgfpointxy{256.029}{302.5}}
    \pgfpatharcaxes{150.0}{180.0}{\pgfpointxy{45.0}{0.0}}{\pgfpointxy{0.0}{45.0}}
    \pgfusepath{stroke}
  \end{pgfscope}
  \begin{pgfscope}
    \pgfpathmoveto{\pgfpointxy{250.0}{280.0}}
    \pgfpatharcaxes{0.0}{30.0}{\pgfpointxy{45.0}{0.0}}{\pgfpointxy{0.0}{45.0}}
    \pgfusepath{stroke}
  \end{pgfscope}
  \begin{pgfscope}
    \pgfpathmoveto{\pgfpointxy{205.33}{302.821}}
    \pgfpatharcaxes{176.565}{206.565}{\pgfpointxy{45.0}{0.0}}{\pgfpointxy{0.0}{45.0}}
    \pgfusepath{stroke}
  \end{pgfscope}
  \begin{pgfscope}
    \pgfpathmoveto{\pgfpointxy{210.0}{280.0}}
    \pgfpatharcaxes{26.5651}{56.5651}{\pgfpointxy{45.0}{0.0}}{\pgfpointxy{0.0}{45.0}}
    \pgfusepath{stroke}
  \end{pgfscope}
  \begin{pgfscope}
    \pgfsetfillcolor{white}
    \pgfpathellipse{\pgfpointxy{250.0}{200.0}}{\pgfpointxy{8.0}{0.0}}{\pgfpointxy{0.0}{8.0}}
    \pgfusepath{fill,stroke}
  \end{pgfscope}
  \begin{pgfscope}
    \pgfsetfillcolor{white}
    \pgfpathellipse{\pgfpointxy{250.0}{400.0}}{\pgfpointxy{8.0}{0.0}}{\pgfpointxy{0.0}{8.0}}
    \pgfusepath{fill,stroke}
  \end{pgfscope}
  \begin{pgfscope}
    \pgfsetfillcolor{white}
    \pgfpathellipse{\pgfpointxy{150.0}{400.0}}{\pgfpointxy{8.0}{0.0}}{\pgfpointxy{0.0}{8.0}}
    \pgfusepath{fill,stroke}
  \end{pgfscope}
  \begin{pgfscope}
    \pgfsetfillcolor{white}
    \pgfpathellipse{\pgfpointxy{350.0}{400.0}}{\pgfpointxy{8.0}{0.0}}{\pgfpointxy{0.0}{8.0}}
    \pgfusepath{fill,stroke}
  \end{pgfscope}
  \begin{pgfscope}
    \pgfsetfillcolor{white}
    \pgfpathellipse{\pgfpointxy{650.0}{200.0}}{\pgfpointxy{8.0}{0.0}}{\pgfpointxy{0.0}{8.0}}
    \pgfusepath{fill,stroke}
  \end{pgfscope}
  \begin{pgfscope}
    \pgfsetfillcolor{white}
    \pgfpathellipse{\pgfpointxy{600.0}{400.0}}{\pgfpointxy{8.0}{0.0}}{\pgfpointxy{0.0}{8.0}}
    \pgfusepath{fill,stroke}
  \end{pgfscope}
  \begin{pgfscope}
    \pgfsetfillcolor{white}
    \pgfpathellipse{\pgfpointxy{500.0}{400.0}}{\pgfpointxy{8.0}{0.0}}{\pgfpointxy{0.0}{8.0}}
    \pgfusepath{fill,stroke}
  \end{pgfscope}
  \begin{pgfscope}
    \pgfsetfillcolor{white}
    \pgfpathellipse{\pgfpointxy{800.0}{400.0}}{\pgfpointxy{8.0}{0.0}}{\pgfpointxy{0.0}{8.0}}
    \pgfusepath{fill,stroke}
  \end{pgfscope}
  \begin{pgfscope}
    \pgfsetfillcolor{white}
    \pgfpathellipse{\pgfpointxy{700.0}{400.0}}{\pgfpointxy{8.0}{0.0}}{\pgfpointxy{0.0}{8.0}}
    \pgfusepath{fill,stroke}
  \end{pgfscope}
  \begin{pgfscope}
    \pgfsetfillcolor{white}
    \pgfpathellipse{\pgfpointxy{1000.0}{200.0}}{\pgfpointxy{8.0}{0.0}}{\pgfpointxy{0.0}{8.0}}
    \pgfusepath{fill,stroke}
  \end{pgfscope}
  \begin{pgfscope}
    \pgfsetfillcolor{white}
    \pgfpathellipse{\pgfpointxy{1050.0}{400.0}}{\pgfpointxy{8.0}{0.0}}{\pgfpointxy{0.0}{8.0}}
    \pgfusepath{fill,stroke}
  \end{pgfscope}
  \begin{pgfscope}
    \pgfsetfillcolor{white}
    \pgfpathellipse{\pgfpointxy{950.0}{400.0}}{\pgfpointxy{8.0}{0.0}}{\pgfpointxy{0.0}{8.0}}
    \pgfusepath{fill,stroke}
  \end{pgfscope}
  \begin{pgfscope}
    \pgfsetfillcolor{white}
    \pgfpathellipse{\pgfpointxy{100.0}{200.0}}{\pgfpointxy{8.0}{0.0}}{\pgfpointxy{0.0}{8.0}}
    \pgfusepath{fill,stroke}
  \end{pgfscope}
  \pgftext[bottom,at={\pgfpointxy{150.0}{420.0}}]{$a_{j_1}$}
  \pgftext[bottom,at={\pgfpointxy{500.0}{420.0}}]{$a_{j_2}$}
  \pgftext[bottom,at={\pgfpointxy{950.0}{420.0}}]{$a_{j_3}$}
  \pgftext[top,at={\pgfpointxy{250.0}{175.0}}]{$h_1(1) = \dots = h_3(1)$}
  \pgftext[top,at={\pgfpointxy{650.0}{175.0}}]{$h_4(1) = \dots = h_7(1)$}
  \pgftext[top,at={\pgfpointxy{1000.0}{175.0}}]{$h_8(1) = h_9(1)$}
  \pgftext[bottom,at={\pgfpointxy{250.0}{487.0}}]{$B_1$}
  \pgftext[bottom,at={\pgfpointxy{650.0}{487.0}}]{$B_2$}
  \pgftext[bottom,at={\pgfpointxy{1000.0}{487.0}}]{$B_3$}
  \pgftext[right,at={\pgfpointxy{63.0}{425.0}}]{$A$}
  \pgftext[right,at={\pgfpointxy{63.0}{300.0}}]{$f^*$}
  \pgftext[right,at={\pgfpointxy{63.0}{200.0}}]{$\omega$}
  \pgftext[left,at={\pgfpointxy{1112.0}{200.0}}]{$\dots$}
\end{pgfpicture}
  \caption{The construction of $f^* : A_\ell \to \omega$}
  \label{rpemklei.fig.snopovi}
\end{figure}

  Claim 1. $\pi$ is injective.

  Proof. Assume that $\pi(g_1) = \pi(g_2) = f^*$ where $f^* : A_\ell \hookrightarrow \omega$ for some
  $A_\ell = \{a_{j_1} < \dots < a_{j_m}\} \subseteq A$ with $j_1 = 1$. Then $g_1$ and $g_2$ are
  coalgebraic homomorphisms between the unordered coalgebras $(A, \alpha)$ and $\calE(\omega)$
  constructed for the $\Set$-monad $(E, \delta, \epsilon)$, Example~\ref{rpemklei.ex.mon-act}.
  Define a mapping $f : A \to \omega$ so that
  \begin{align*}
    f^*(a_{1})   &= f(a_{j_1}) = f(a_{j_1+1}) = \dots = f(a_{j_2 - 1}), && \text{[note: $j_1 = 1$]}\\
    f^*(a_{j_2}) &= f(a_{j_2}) = f(a_{j_2 + 1}) = \dots = f(a_{j_3 - 1}), \\
                 &\;\vdots \\
    f^*(a_{j_m}) &= f(a_{j_m}) = f(a_{j_m + 1}) = \dots = f(a_{s}).
  \end{align*}
  Then the definition of $\pi$ implies that $\epsilon_\omega \cdot g_1 = f$ and $\epsilon_\omega \cdot g_2 = f$.
  Since $\calE(X)$ is a cofree $E$-coalgebra, $g_1 = g_2$ by Lemma~\ref{rpemklei.lem.uniq-cofree}.
  This concludes the proof of Claim~1.

  \medskip

  Define $\gamma : \pi(R) \to k$ by $\gamma(\pi(f)) = \chi(f)$ so that
  $
    \chi(R) = \gamma(\pi(R))
  $
  and then define $\gamma_i : S_i \to k$, $1 \le i \le n$, by
  $$
    \gamma_i(h) = \begin{cases}
      \gamma(h), & h \in \pi(R) \cap S_i,\\
      0, & \text{otherwise}.
    \end{cases}
  $$
  Let us construct
  $\gamma'_i : S_i \to k$ and $w_i \in \hom_{\ChEmb}(\omega, \omega)$, $i \in \{1, \dots, n\}$,
  inductively as follows. First, put $\gamma'_n = \gamma_n$.
  Given a coloring $\gamma'_i : S_i \to k$, construct $w_i$
  by the Infinite Ramsey Theorem (Theorem~\ref{rpemklei.thm.IRT}):
  since $\omega \longrightarrow (\omega)^{A_i}_{k}$, there is a $w_i \in \hom_{\ChEmb}(\omega, \omega)$
  such that
  $
    |\gamma'_i(w_i \cdot S_i)| \le 1
  $.
  Finally, given $w_i \in \hom_{\ChEmb}(\omega, \omega)$ define $\gamma'_{i-1} : S_{i-1} \to k$ by
  $
    \gamma'_{i-1}(f) = \gamma_{i-1}(w_n \cdot \ldots \cdot w_i \cdot f)
  $.
  Now, put $u = w_n \cdot \ldots \cdot w_1 \in \hom_{\ChEmb}(\omega,\omega)$ and let us show that
  $
    |\chi(\hat E(u) \cdot R)| \le n
  $.

  \medskip

  Claim 2. $\pi(\hat E(u) \cdot R) = u \cdot \pi(R) \subseteq \pi(R)$.

  Proof. Since $\hat E(u) \cdot R \subseteq R$ it immediately follows that $\pi(\hat E(u) \cdot R) \subseteq \pi(R)$.
  To see that $\pi(\hat E(u) \cdot R) = u \cdot \pi(R)$ take any $f \in R$, let $g = \hat E(u) \cdot f$ and let us show
  that $g^* = u \cdot f^*$. Following the definition of $f^*$ let $f(a_i) = h_i$, $1 \le i \le s$.
  Then $g(a_i) = (\hat E(u) \cdot f)(a_i) = \hat E(u)(f(a_i)) = \hat E(u)(h_i) = u \cdot h_i$.
  As above, we easily conclude that $h_1(1) \le h_2(1) \le \dots \le h_s(1)$ because
  $f$ is an embedding and $\omega^M$ is ordered lexicographically. Since
  $u : \omega \to \omega$ is an embedding, $u \cdot h_1(1) \le u \cdot h_2(1) \le \dots \le u \cdot h_s(1)$.
  Moreover, $h_i(1) = h_j(1)$ iff $u \cdot h_i(1) = u \cdot h_j(1)$ for all $1 \le i, j \le s$. The definition of $f^*$
  then immediately gives us that $g^* = u \cdot f^*$. This concludes the proof of the claim.

  \medskip

  Therefore, $\chi(\hat E(u) \cdot R) = \gamma(\pi(\hat E(u) \cdot R)) = \gamma(u \cdot \pi(R))$.
  Since $\pi(R) \subseteq \bigcup_{i=1}^n S_i$,
  $$
    |\gamma(u \cdot \pi(R))|
    \le |\gamma(u \cdot \bigcup_{i=1}^n S_i)|
    = |\bigcup_{i=1}^n \gamma(u \cdot S_i)|
    \le \sum_{i=1}^n |\gamma(u \cdot S_i)|.
  $$
  Fix an $i \in \{1, \dots, n\}$. Clearly, $u \cdot S_i \subseteq S_i$
  and $w_i \cdot \ldots \cdot w_1 \cdot S_i  \subseteq w_i \cdot S_i$, whence
  $$
    |\gamma(u \cdot S_i)|
    = |\gamma_i(w_n \cdot \ldots \cdot w_1 \cdot S_i)|
    = |\gamma'_i(w_i \cdot \ldots \cdot w_1 \cdot S_i)|
    \le |\gamma'_i(w_i \cdot S_i)| \le 1.
  $$
  Putting it all together, we finally get
  $
    |\chi(\hat E(u) \cdot R)| = |\gamma(u \cdot \pi(R))| \le \sum_{i=1}^n |\gamma(u \cdot S_i)| \le n = 2^{s-1}
  $.
  This completes the proof.
\end{proof}

\begin{corol}
  Let $M$ be an arbitrary monoid (finite or infinite).
  Every finite $M$-set has a finite big Ramsey degree in $\calE(\omega)$, the cofree $M$-set on $\omega$ generators.
  More precisely, for every finite $M$-set $\calA$ with $n$ elements we have
  $T(\calA, \calE(\omega)) \le n! \cdot 2^{n-1}$ in $\SetEmb(M)$.
  
  In particular, for every group $G$ every finite $G$-set has a finite big Ramsey degree in the cofree $G$-set on $\omega$ generators;
  and for every unary algebraic language $\Omega$ every finite $\Omega$-algebra
  has a finite big Ramsey degree in the cofree $\Omega$-algebra on $\omega$ generators.
\end{corol}
\begin{proof}
  As a matter of notational convenience put $\CCC = \OsetEmb(M)$ and $\BB = \SetEmb(M)$.
  Let $U : \CCC^\op \to \BB^\op$ be the forgetful functor that forgets the order.
  To see that $U$ has unique restrictions we can repeat the argument from the proof of
  Corollary~\ref{rpemklei.cor.fsrd-G-set}.

  Let $\calA = (A, \Omega^A) \in \Ob(\BB^\fin)$ be an arbitrary finite $M$-set.
  As we have seen in Theorem~\ref{rpemklei.thm.main-G-set-order},
  there exists a linear ordering $\hat\calE(\omega)$ of $\calE(\omega)$
  such that every finite ordered $M$-set with $n$ elements has a big dual Ramsey degree in $\hat\calE(\omega)$
  which does not exceed~$2^{n-1}$. Therefore, for every linear ordering $<$ of $A$ we have that
  $
    T_\CCC((A, \alpha, <), \hat\calE(\omega)) \le 2^{n-1}
  $.
  Since $U^{-1}(\calA)$ is finite (because there are only finitely many linear orders on a finite set)
  Theorem~\ref{rpemklei---big.v.small} tells us that
  $
    T_{\BB}(\calA, \calE(\omega)) \le \sum_{\calA^* \in U^{-1}(\calA)} T_{\CCC}(\calA^*, \hat\calE(\omega)) \le n! \cdot 2^{n-1}
  $,
  where $n = |A|$. This completes the proof.
\end{proof}


\end{document}